\documentclass[a4paper]{article}
\makeatletter
\providecommand{\iftx@stxtwo}{\iffalse}
\providecommand{\iftx@ebgm}{\iffalse}
\providecommand{\iftx@ut}{\iffalse}
\providecommand{\iftx@nc}{\iftrue}
\providecommand{\iftx@fourier}{\iffalse}
\makeatother

\usepackage{xparse}
\usepackage[intlimits]{amsmath}
\usepackage[lsym=ptm,cal=dutch,scr=true]{mwmfonts}
\usepackage{amssymb,amsfonts,amsthm,mathtools}
\usepackage{tikz}
\usetikzlibrary{decorations.markings,decorations.pathreplacing,positioning,cd,fit,shapes.geometric}
\tikzset{radial vector/.style={blue, -latex}}
\tikzset{radial vector 2/.style={red, -latex}}
\tikzset{midarrow/.style={postaction=decorate,decoration={markings,mark={at position #1 with {\arrow{latex}}}}}}
\tikzset{invmidarrow/.style={postaction=decorate,decoration={markings,mark={between positions 0 and 1 step #1 with {\arrow{latex reversed}}}}}}
\usepackage{cancel}
\usepackage[top=3cm, bottom=3.5cm, outer=4cm, inner=4cm, marginparsep=0.7cm, marginparwidth=1.5cm]{geometry}
\usepackage{cancel,hyperref}
\usepackage{newtxtextM,color,thmtools}

\newcommand{\nipar}{\par\noindent}
\makeatletter

\let\eatup\@gobble
\newcommand*\MapsTo{%
  \@ifstar{\xrightarrow[\raisebox{0.25 em}{\smash{\ensuremath{\sim}}}]{}}{\xrightarrow{\raisebox{-0.25 em}{\smash{\ensuremath{\sim}}}}}%
}
\newcommand{\fadecol}[1]{\textcolor{gray}{#1}}
\newif\ifmidv\midvtrue
\newcommand{\midv}{
    \newcommand{\fade}[1]{\fadecol{##1}}
    \newcommand{\fadees}{\fade}
    \newcommand{\fadecap}[3]{\fadees{##1}{##2}{\MakeUppercase##3}}
    \newcommand{\lac}[1]{\@ifstar{<>##1}{<##1>}}
    \newcommand{\slac}[1]{<##1> }
    \newcommand{\xlac}[1]{\textcolor{red}{<}##1\textcolor{red}{>}\@ifstar{\ }{}}
    \newcommand{\xsqb}{\sqb}
    \midvtrue
}
\newcommand{\finv}{
    \newcommand{\fade}[1]{}
    \newcommand{\fadees}{\fade}
    \newcommand{\fadecap}[3]{\fadees{##1}{##2}{\MakeUppercase##3}}
    \newcommand{\lac}{\@ifstar{\expandafter\@gobble\@gobble}{}}
    \newcommand{\slac}{\lac}
    \newcommand{\xlac}{\@ifstar{\lac}{\lac}}
    \newcommand{\xsqb}{}
    \finvtrue
}
\makeatother
\newcommand{\opn}{\operatorname}
\catcode`\à=\active
\defà{\`a}
\catcode`\è=\active
\defè{\`e}
\catcode`\ì=\active
\defì{\`\i}
\catcode`\ò=\active
\defò{\`o}
\catcode`\ù=\active
\defù{\`u}
\catcode`\á=\active
\defá{\'a}
\catcode`\é=\active
\defé{\'e}
\catcode`\í=\active
\defí{\'\i}
\catcode`\ó=\active
\defó{\'o}
\catcode`\ú=\active
\defú{\'u}
\catcode`\ã=\active
\defã{\txtilde{a}}
\catcode`\ẽ=\active
\defẽ{\txtilde{e}}
\catcode`\ĩ=\active
\defĩ{\txtilde{\i}}
\catcode`\õ=\active
\defõ{\txtilde{o}}
\catcode`\ũ=\active
\defũ{\txtilde{u}}
\catcode`\â=\active
\defâ{\^a}
\catcode`\ê=\active
\defê{\^e}
\catcode`\î=\active
\defî{\^\i}
\catcode`\ô=\active
\defô{\^o}
\catcode`\û=\active
\defû{\^u}
\catcode`\č=\active
\defč{\v c}
\catcode`\š=\active
\defš{\v s}
\catcode`\Č=\active
\defČ{\v C}
\catcode`\Š=\active
\defŠ{\v S}
\catcode`\ß=\active
\defß{\sz}
\catcode`\ä=\active
\defä{\"a}
\catcode`\ë=\active
\defë{\"e}
\catcode`\ï=\active
\defï{\"\i}
\catcode`\ö=\active
\defö{\"o}
\catcode`\ü=\active
\defü{\"u}
\catcode`\Ä=\active
\defÄ{\"A}
\catcode`\Ë=\active
\defË{\"E}
\catcode`\Ï=\active
\defÏ{\"I}
\catcode`\Ö=\active
\defÖ{\"O}
\catcode`\Ü=\active
\defÜ{\"U}
\catcode`\ā=\active
\defā{\=a}
\catcode`\ē=\active
\defē{\=e}
\catcode`\ī=\active
\defī{\=\i}
\catcode`\ō=\active
\defō{\=o}
\catcode`ū=\active
\defū{\=u}
\catcode`\Ā=\active
\defĀ{\=A}
\catcode`\Ē=\active
\defĒ{\=E}
\catcode`\Ī=\active
\defĪ{\=I}
\catcode`\Ō=\active
\defŌ{\=O}
\catcode`\Ū=\active
\defŪ{\=U}
\catcode`\ç=\active
\defç{\c c}
\catcode`\ğ=\active
\defğ{\u g}
\catcode`\ı=\active
\defı{\i}
\catcode`\ş=\active
\defş{\c s}
\catcode`\Ç=\active
\defÇ{\c C}
\catcode`\Ğ=\active
\defĞ{\u G}
\catcode`\İ=\active
\defİ{\. I}
\catcode`\Ş=\active
\defŞ{\c S}
\catcode`\ś=\active
\defś{\'s}
\catcode`\Ś=\active
\defŚ{\'S}
\catcode`\ć=\active
\defć{\'c}
\catcode`\Ć=\active
\defĆ{\'C}
\catcode`\ǵ=\active
\defǵ{\'g}
\catcode`\Ǵ=\active
\defǴ{\'G}
\catcode`\ŝ=\active
\defŝ{\^s}
\catcode`\Ŝ=\active
\defŜ{\^S}
\catcode`\Ṣ=\active
\defṢ{\d S}
\catcode`\ṣ=\active
\defṣ{\d s}
\catcode`\Ǒ=\active
\defǑ{\v O}
\catcode`\ǒ=\active
\defǒ{\v o}
\catcode`\Ã=\active
\defÃ{\txtilde A}
\catcode`\Á=\active
\defÁ{\' A}
\catcode`\À=\active
\defÀ{\` A}
\catcode`\Ẽ=\active
\defẼ{\txtilde E}
\catcode`\É=\active
\defÉ{\' E}
\catcode`\Ê=\active
\defÊ{\^ E}
\catcode`\È=\active
\defÈ{\`E}
\catcode`\Ì=\active
\defÌ{\` I}
\catcode`\Õ=\active
\defÕ{\txtilde O}
\catcode`\Ò=\active
\defÒ{\` O}
\catcode`\Ó=\active
\defÓ{\' O}
\catcode`\Ô=\active
\defÔ{\^O}
\catcode`\Ù=\active
\defÙ{\`U}
\catcode`\Ẓ=\active
\defẒ{\d Z}
\catcode`\ẓ=\active
\defẓ{\d z}
\catcode`\Ȧ=\active
\defȦ{\.A}
\catcode`\ȧ=\active
\defȧ{\.a}
\catcode`\Ě=\active
\defĚ{\v E}
\catcode`\ě=\active
\defě{\v e}
\newcommand{\grad}{\nabla}
\newcommand{\pd}{\partial}
\newcommand{\sect}{\section}
\newcommand{\ssect}{\subsection}

\newcommand{\ben}{\begin{enumerate}}
\newcommand{\een}{\end{enumerate}}
\newcommand{\bde}{\begin{description}}
\newcommand{\ede}{\end{description}}
\newcommand{\bi}{\begin{itemize}}
\newcommand{\ei}{\end{itemize}}
\newcommand{\hsp}{\hspace}
\newcommand{\vsp}{\vspace}
\newcommand{\fn}{\footnote}

\newcommand{\fr}{\dfrac}
\newcommand{\Der}[2][]{\fr{{\diff #1}}{{\diff #2}}}

\newcommand{\e}[1][2]{\mskip-#1\thinmuskip}
\newcommand{\diff}{\mathrm{d}\mskip.25\thinmuskip}

\newcommand{\autocnt}{\@ifstar{\@starcnttrue\@autocnt}{\@starcntfalse\@autocnt}}
\newcommand{\refprep}{di}
\newcommand{\setprep}{\renewcommand\refprep}
\newcommand{\refpunct}{:}
\makeatletter
\newif\ifinvd@lims\invd@limsfalse
\newif\iflr\lrfalse
\newif\if@labeling\@labelingfalse
\newif\ifplaceqed\placeqedfalse
\newif\if@starcnt\@starcnttrue
\newif\if@namespace\@namespacefalse
    \newcommand{\mw@contname}{continues~from}
    \newcommand{\mw@Contname}{Continues~from}
    \newcommand{\mw@definame}{Definition}
    \newcommand{\mw@teorname}{Theorem}
    \newcommand{\mw@esname}{Example}
    \newcommand{\mw@esename}{Exercise}
    \newcommand{\mw@ptname}{Point}
    \newcommand{\mw@riassname}{Sum-up}
    \newcommand{\mw@proponame}{Proposition}
    \newcommand{\mw@ossname}{Remark}
    \newcommand{\mw@corname}{Corollary}
    \newcommand{\mw@esecasaname}{Exercise~done~by~me}
    \newcommand{\mw@lemmaname}{Lemma}
    \newcommand{\mw@fattoname}{Fact}
    \newcommand{\mw@proofstring}{Proof.}
    \newcommand{\mw@lessname}{Lessons~by}
    \newcommand{\mw@esercname}{Exercise~lessons~by}
    \newcommand{\mw@notename}{Notes~by}
    \newcommand{\mw@acyearname}{Academic~Year}
    \newcommand{\mw@alsoname}{And~also~by}
    \setprep{of}

\newcommand{\autocnts}{\@ifstar{\@starcnttrue\@autocnts}{\@starcntfalse\@autocnts}}
\newcommand{\@autocnts}{\@autocnt{defi}\@autocnt{oss}\@autocnt{propo}\@autocnt{lemma}\@autocnt{es}\@autocnt{teor}\@autocnt{ese}\@autocnt{cor}\@autocnt{fatto}}
\newcommand{\xhapter}[1]{\@ifstar{\chapter*{#1}\autocnts}{\chapter{#1}\autocnts}}
\newcommand{\xect}[1]{\@ifstar{\section*{#1}\autocnts}{\section{#1}\autocnts}}
\newcommand{\xsect}[1]{\@ifstar{\subsection*{#1}\autocnts}{\subsection{#1}\autocnts}}
\newcommand{\@autocnt}[1]
    {\if@starcnt\setcounter{#1}{0}\fi
    \ifnum\c@chapter>0
        \ifnum\c@section>0
            \ifnum\c@subsection>0
                \expandafter\renewcommand\csname the#1\endcsname{\thesubsection.\arabic{#1}%
                \numberwithin{#1}{subsection}}%
            \else
                \expandafter\renewcommand\csname the#1\endcsname{\thesection.\arabic{#1}%
                \numberwithin{#1}{section}}%
            \fi
        \else
            \expandafter\renewcommand\csname the#1\endcsname{\thechapter.\arabic{#1}%
                \numberwithin{#1}{chapter}}%
         \fi
     \else
         \expandafter\renewcommand\csname the#1\endcsname{\arabic{#1}}%
     \fi
     }
\newcommand{\lonecnt}{\@ifstar{\@starcnttrue\@lonecnt}{\@starcntfalse\@lonecnt}}
\newcommand{\@lonecnt}[1]{\if@starcnt\setcounter{#1}{0}\fi
\expandafter\renewcommand\csname the#1\endcsname{\arabic{#1}}}
\newcommand{\nocnt}{\@ifstar{\@starcnttrue\@nocnt}{\@starcntfalse\@nocnt}}
\newcommand{\@nocnt}[1]{\if@starcnt\setcounter{#1}{0}\fi
\expandafter\def\csname the#1\endcsname{\hsp{\mw@No}}}
\newlength\mw@No
\newcommand{\setNo}[1]{\setlength{\mw@No}{#1}}
\newcommand{\resetNo}[1]{\setlength{\mw@No}{-3pt}}
\newcommand{\chaptercnt}{\@ifstar{\@starcnttrue\@chaptercnt}{\@starcntfalse\@chaptercnt}}
\newcommand{\@chaptercnt}[1]{\if@starcnt\setcounter{#1}{0}\fi
\expandafter\renewcommand\csname the#1\endcsname{\thechapter.\arabic{#1}}\numberwithin{#1}{chapter}}
\newcommand{\sectcnt}{\@ifstar{\@starcnttrue\@sectcnt}{\@starcntfalse\@sectcnt}}
\newcommand{\@sectcnt}[1]{\if@starcnt\setcounter{#1}{0}\fi
\expandafter\renewcommand\csname the#1\endcsname{\thesection.\arabic{#1}}\numberwithin{#1}{section}}
\newcommand{\ssectcnt}{\@ifstar{\@starcnttrue\@ssectcnt}{\@starcntfalse\@ssectcnt}}
\newcommand{\@ssectcnt}[1]{\if@starcnt\setcounter{#1}{0}\fi
\if@starcnt\setcounter{#1}{0}\fi
\expandafter\renewcommand\csname the#1\endcsname{\thesubsection.\arabic{#1}}\numberwithin{#1}{subsection}}
\newcommand{\xtheorsetlist}[1][8]{\renewcommand\thmt@listnumwidth{#1em}}
\ExplSyntaxOn
\NewDocumentCommand{\definethm}{smmO{a}}
    {\newtheorem{#2}{#3}
        \IfBooleanTF{#1}
            {\setcnttype*{#2}{#4}}
            {\setcnttype{#2}{#4}}
        \expandafter\newcommand\csname mw@#2name\endcsname{#3}
    }
\NewDocumentCommand{\xbegin}{moood()}
    {\IfValueTF{#5}
        {\str_case:nnF{#2}{
        {\ }{\gdef\thisthmname{\protect\nameref*{#3}}}
        {n}{\let\thisthmname\relax}
        }{\gdef\thisthmname{#2}}
        \begin{#1}[continued=#3,contfrom=#5,label=#4]}
        {\IfValueTF{#4}
            {\str_if_eq:nnTF{#2}{\ }
                {\begin{#1}[{name=[#4]},label={#3}]}
                {\begin{#1}[{name=[#4]#2},label={#3}]}}
            {\IfValueTF{#3}
                {\str_if_eq:nnTF{#2}{\ }
                 {\begin{#1}[label={#3}]}
                 {\begin{#1}[name={#2},label={#3}]}}
                {\IfValueTF{#2}
                    {\begin{#1}[name={#2}]}
                    {\begin{#1}}
                }
            }
        }
    }
\newcommand{\xend}[1]{\end{#1}\noindent\ignorespaces}
\newcommand{\setcnttype}{\@ifstar{\@starcnttrue\@setcnttype}{\@starcntfalse\@setcnttype}}
\newcommand{\@setcnttype}[2]
    {\str_case:nnF{#2}{
        {a}{\@autocnt{#1}}
        {n}{\@lonecnt{#1}}
        {@@}{\@nocnt{#1}}
        {c}{\@chaptercnt{#1}}
        {s}{\@sectcnt{#1}}
        {x}{\@ssectcnt{#1}}
        {nc}{}
        }
    {\msg_fatal:nnn{mworksx}{cnt-type-undef}{#2}}
}
\NewDocumentEnvironment{qeddim}{o}%
    {\global\placeqedfalse\noindent
    \IfValueTF{#1}
        {#1 \\ }
        {}
    \mw@proofstring\ \par}
    {\ifplaceqed\else\hfill\qedsym\linebreak\fi}
\NewDocumentEnvironment{qeddim*}{o}%
    {\global\placeqedfalse\noindent
    \mw@proofstring\IfValueTF{#1}{\ (#1)}{}\par}
    {\ifplaceqed\else\hfill\qedsym\linebreak\fi}
\DeclareSymbolFont{lasy}{U}{lasy}{m}{n}
\SetSymbolFont{lasy}{bold}{U}{lasy}{b}{n}
\DeclareMathSymbol\Diamondx{\mathord}{lasy}{"33}
\newcommand{\qedsym}{$\Diamondx$}

\newcommand{\placeqed}{\@ifstar{\hfill\qedsym\global\placeqedtrue}{\tag*{\qedsym}\global\placeqedtrue}}

\NewDocumentCommand{\@ag}{O{0.85}oO{}}
    {\str_if_eq:nnTF{#1}{n}
            {\IfValueTF{#2}
                {#2{\scalebox{0.85}[0.85]{$#3\alpha$}}}
                {\errmessage{Too few arguments}\errhelp{If you specify `n' as the first argument, the second optional argument must be specified, otherwise use a number as the first argument or specify no optional argument.}}
            }{\str_if_eq:nnTF{#1}{e}
                {\IfValueTF{#2}
                    {#2{\scalebox{0.65}[0.65]{$#3\alpha$}}}
                    {\scalebox{0.65}[0.65]{$#3\alpha$}}
            }{\str_if_eq:nnTF{#1}{b}
                    {\IfValueTF{#2}
                        {#2{#3\alpha}}
                        {{#3\alpha}}
             }{\IfValueTF{#2}
                        {#2{\scalebox{#1}[#1]{$#3\alpha$}}}
                        {\scalebox{#1}[#1]{$#3\alpha$}}
    }}}}

\newcommand{\ag}{
    \@ifgrave
        {
        \@ifstar
            {
            \@ag[e][\mathclap]
            }
            {
            \@ag[e]
            }
        }
        {    \@ifstar
             {
             \@ag[n][\mathclap]
             }
             {
             \@ag[n][]
             }
         }
     }
     
\newcommand{\xag}{\@ag}

\NewDocumentCommand{\@bg}{O{0.88}oO{}}
    {\str_if_eq:nnTF{#1}{n}
            {\IfValueTF{#2}
                {#2{\scalebox{0.88}[0.88]{$#3\beta$}}}
                {\errmessage{Too few arguments}\errhelp{If you specify `n' as the first argument, the second optional argument must be specified, otherwise use a number as the first argument or specify no optional argument.}}
            }{\str_if_eq:nnTF{#1}{e}
                {\IfValueTF{#2}
                    {#2{\scalebox{0.65}[0.65]{$#3\beta$}}}
                    {\scalebox{0.65}[0.65]{$#3\beta$}}
                }{\str_if_eq:nnTF{#1}{b}
                    {\IfValueTF{#2}
                        {#2{#3\beta}}
                        {{#3\beta}}
                    }{\IfValueTF{#2}
                        {#2{\scalebox{#1}[#1]{$#3\beta$}}}
                        {\scalebox{#1}[#1]{$#3\beta$}}
                    }
                }
            }
    }

\newcommand{\bg}{
    \@ifgrave
        {
        \@ifstar
            {
            \@bg[e][\mathclap]
            }
            {
            \@bg[e]
            }
        }
        {    \@ifstar
             {
             \@bg[n][\mathclap]
             }
             {
             \@bg
             }
         }
     }
     
\newcommand{\xbg}{\@bg}

\NewDocumentCommand{\@gg}{O{0.85}oO{}}
{\str_if_eq:nnTF{#1}{n}
        {\IfValueTF{#2}
            {#2{\scalebox{0.85}[0.85]{$#3\gamma$}}}
            {\errmessage{Too few arguments}\errhelp{If you specify `n' as the first argument, the second optional argument must be specified, otherwise use a number as the first argument or specify no optional argument.}}
        }{\str_if_eq:nnTF{#1}{e}
            {\IfValueTF{#2}
                {#2{\scalebox{0.65}[0.65]{$#3\gamma$}}}
                {\scalebox{0.65}[0.65]{$#3\gamma$}}
            }{\str_if_eq:nnTF{#1}{b}
                {\IfValueTF{#2}
                    {#2{#3\gamma}}
                    {{#3\gamma}}
               }{\str_if_eq:nnTF{#1}{t}
                     {\IfValueTF{#2}
                         {#2{\scalebox{0.88}[0.88]{$#3\gamma$}}}
                         {\scalebox{0.88}[0.88]{$#3\gamma$}}
                     }{\IfValueTF{#2}
                        {#2{\scalebox{#1}[#1]{$#3\gamma$}}}
                        {\scalebox{#1}[#1]{$#3\gamma$}}
}}}}}

\def\gg{
    \@ifgrave
        {
        \@ifstar
            {
            \@gg[e][\mathclap]
            }
            {
            \@gg[e]
            }
        }
        {    \@ifstar
             {
             \@gg[n][\mathclap]
             }
             {
             \@gg
             }
         }
     }
     
\newcommand{\xgg}{\@gg}

\NewDocumentCommand{\@dg}{O{0.88}oO{}}
{\str_if_eq:nnTF{#1}{n}
        {\IfValueTF{#2}
            {#2{\scalebox{0.88}[0.88]{$#3\delta$}}}
            {\errmessage{Too few arguments}\errhelp{If you specify `n' as the first argument, the second optional argument must be specified, otherwise use a number as the first argument or specify no optional argument.}}
        }{\str_if_eq:nnTF{#1}{e}
            {\IfValueTF{#2}
                {#2{\scalebox{0.65}[0.65]{$#3\delta$}}}
                {\scalebox{0.65}[0.65]{$#3\delta$}}
        }{\str_if_eq:nnTF{#1}{b}
                {\IfValueTF{#2}
                    {#2{#3\delta}}
                    {{#3\delta}}
         }{\IfValueTF{#2}
                    {#2{\scalebox{#1}[#1]{$#3\delta$}}}
                    {\scalebox{#1}[#1]{$#3\delta$}}
}}}}

\newcommand{\dg}{
    \@ifgrave
        {
        \@ifstar
            {
            \@dg[e][\mathclap]
            }
            {
            \@dg[e]
            }
        }
        {    \@ifstar
             {
             \@dg[n][\mathclap]
             }
             {
             \@dg
             }
         }
     }
     
\newcommand{\xdg}{\@dg}

\NewDocumentCommand{\@eg}{O{0.85}oO{}}
{\str_if_eq:nnTF{#1}{n}
        {\IfValueTF{#2}
            {#2{\scalebox{0.85}[0.85]{$#3\varepsilon$}}}
            {\errmessage{Too few arguments}\errhelp{If you specify `n' as the first argument, the second optional argument must be specified, otherwise use a number as the first argument or specify no optional argument.}}
        }{\str_if_eq:nnTF{#1}{e}
            {\IfValueTF{#2}
                {#2{\scalebox{0.65}[0.65]{$#3\varepsilon$}}}
                {\scalebox{0.65}[0.65]{$#3\varepsilon$}}
        }{\str_if_eq:nnTF{#1}{b}
                {\IfValueTF{#2}
                    {#2{#3\varepsilon}}
                    {{#3\varepsilon}}
         }{\IfValueTF{#2}
                    {#2{\scalebox{#1}[#1]{$#3\varepsilon$}}}
                    {\scalebox{#1}[#1]{$#3\varepsilon$}}
}}}}

\newcommand{\eg}{
    \@ifgrave
        {
        \@ifstar
            {
            \@eg[e][\mathclap]
            }
            {
            \@eg[e]
            }
        }
        {    \@ifstar
             {
             \@eg[n][\mathclap]
             }
             {
             \@eg
             }
         }
     }
     
\newcommand{\xeg}{\@eg}

\NewDocumentCommand{\@zg}{O{0.88}oO{}}
{\str_if_eq:nnTF{#1}{n}
        {\IfValueTF{#2}
            {#2{\scalebox{0.88}[0.88]{$#3\zeta$}}}
            {\errmessage{Too few arguments}\errhelp{If you specify `n' as the first argument, the second optional argument must be specified, otherwise use a number as the first argument or specify no optional argument.}}
        }{\str_if_eq:nnTF{#1}{e}
            {\IfValueTF{#2}
                {#2{\scalebox{0.65}[0.65]{$#3\zeta$}}}
                {\scalebox{0.65}[0.65]{$#3\zeta$}}
        }{\str_if_eq:nnTF{#1}{b}
                {\IfValueTF{#2}
                    {#2{#3\zeta}}
                    {{#3\zeta}}
         }{\IfValueTF{#2}
                    {#2{\scalebox{#1}[#1]{$#3\zeta$}}}
                    {\scalebox{#1}[#1]{$#3\zeta$}}
}}}}

\newcommand{\zg}{
    \@ifgrave
        {
        \@ifstar
            {
            \@zg[e][\mathclap]
            }
            {
            \@zg[e]
            }
        }
        {    \@ifstar
             {
             \@zg[n][\mathclap]
             }
             {
             \@zg
             }
         }
     }

\newcommand{\xzg}{\@zg}

\NewDocumentCommand{\@hg}{O{0.85}oO{}}
{\str_if_eq:nnTF{#1}{n}
        {\IfValueTF{#2}
            {#2{\scalebox{0.85}[0.85]{$#3\eta$}}}
            {\errmessage{Too few arguments}\errhelp{If you specify `n' as the first argument, the second optional argument must be specified, otherwise use a number as the first argument or specify no optional argument.}}
        }{\str_if_eq:nnTF{#1}{e}
            {\IfValueTF{#2}
                {#2{\scalebox{0.65}[0.65]{$#3\eta$}}}
                {\scalebox{0.65}[0.65]{$#3\eta$}}
        }{\str_if_eq:nnTF{#1}{b}
                {\IfValueTF{#2}
                    {#2{#3\eta}}
                    {{#3\eta}}
         }{\str_if_eq:nnTF{#1}{t}
             {\IfValueTF{#2}
                 {#2{\scalebox{0.88}[0.88]{$#3\eta$}}}
                 {\scalebox{0.88}[0.88]{$#3\eta$}}
                 }{\IfValueTF{#2}
                    {#2{\scalebox{#1}[#1]{$#3\eta$}}}
                    {\scalebox{#1}[#1]{$#3\eta$}}
}}}}}

\def\hg{
    \@ifgrave
        {
        \@ifstar
            {
            \@hg[e][\mathclap]
            }
            {
            \@hg[e]
            }
        }
        {    \@ifstar
             {
             \@hg[n][\mathclap]
             }
             {
             \@hg
             }
         }
     }

\newcommand{\xhg}{\@hg}

\NewDocumentCommand{\@qg}{O{0.95}oO{}}
{\str_if_eq:nnTF{#1}{n}
        {\IfValueTF{#2}
            {#2{\scalebox{0.95}[0.95]{$#3\theta$}}}
            {\errmessage{Too few arguments}\errhelp{If you specify `n' as the first argument, the second optional argument must be specified, otherwise use a number as the first argument or specify no optional argument.}}
        }{\str_if_eq:nnTF{#1}{e}
            {\IfValueTF{#2}
                {#2{\scalebox{0.7}[0.7]{$#3\theta$}}}
                {\scalebox{0.7}[0.7]{$#3\theta$}}
        }{\str_if_eq:nnTF{#1}{b}
                {\IfValueTF{#2}
                    {#2{#3\theta}}
                    {{#3\theta}}
         }{\IfValueTF{#2}
                    {#2{\scalebox{#1}[#1]{$#3\theta$}}}
                    {\scalebox{#1}[#1]{$#3\theta$}}
}}}}

\newcommand{\qg}{
    \@ifgrave
        {
        \@ifstar
            {
            \@qg[e][\mathclap]
            }
            {
            \@qg[e]
            }
        }
        {    \@ifstar
             {
             \@qg[n][\mathclap]
             }
             {
             \@qg
             }
         }
     }

\newcommand{\xqg}{\@qg}

\NewDocumentCommand{\@ig}{O{0.85}oO{}}
{\str_if_eq:nnTF{#1}{n}
        {\IfValueTF{#2}
            {#2{\scalebox{0.85}[0.85]{$#3\iota$}}}
            {\errmessage{Too few arguments}\errhelp{If you specify `n' as the first argument, the second optional argument must be specified, otherwise use a number as the first argument or specify no optional argument.}}
        }{\str_if_eq:nnTF{#1}{e}
            {\IfValueTF{#2}
                {#2{\scalebox{0.65}[0.65]{$#3\iota$}}}
                {\scalebox{0.65}[0.65]{$#3\iota$}}
        }{\str_if_eq:nnTF{#1}{b}
                {\IfValueTF{#2}
                    {#2{#3\iota}}
                    {{#3\iota}}
         }{\IfValueTF{#2}
                    {#2{\scalebox{#1}[#1]{$#3\iota$}}}
                    {\scalebox{#1}[#1]{$#3\iota$}}
}}}}

\newcommand{\ig}{
    \@ifgrave
        {
        \@ifstar
            {
            \@ig[e][\mathclap]
            }
            {
            \@ig[e]
            }
        }
        {    \@ifstar
             {
             \@ig[n][\mathclap]
             }
             {
             \@ig
             }
         }
     }

\newcommand{\xig}{\@ig}

\NewDocumentCommand{\@kg}{O{0.85}oO{}}
{\str_if_eq:nnTF{#1}{n}
        {\IfValueTF{#2}
            {#2{\scalebox{0.85}[0.85]{$#3\kappa$}}}
            {\errmessage{Too few arguments}\errhelp{If you specify `n' as the first argument, the second optional argument must be specified, otherwise use a number as the first argument or specify no optional argument.}}
        }{\str_if_eq:nnTF{#1}{e}
            {\IfValueTF{#2}
                {#2{\scalebox{0.65}[0.65]{$#3\kappa$}}}
                {\scalebox{0.65}[0.65]{$#3\kappa$}}
        }{\str_if_eq:nnTF{#1}{b}
                {\IfValueTF{#2}
                    {#2{#3\kappa}}
                    {{#3\kappa}}
         }{\IfValueTF{#2}
                    {#2{\scalebox{#1}[#1]{$#3\kappa$}}}
                    {\scalebox{#1}[#1]{$#3\kappa$}}
}}}}

\newcommand{\kg}{
    \@ifgrave
        {
        \@ifstar
            {
            \@kg[e][\mathclap]
            }
            {
            \@kg[e]
            }
        }
        {    \@ifstar
             {
             \@kg[n][\mathclap]
             }
             {
             \@kg
             }
         }
     }

\newcommand{\xkg}{\@kg}

\NewDocumentCommand{\@lg}{O{0.88}oO{}}
{\str_if_eq:nnTF{#1}{n}
        {\IfValueTF{#2}
            {#2{\scalebox{0.88}[0.88]{$#3\lambda$}}}
            {\errmessage{Too few arguments}\errhelp{If you specify `n' as the first argument, the second optional argument must be specified, otherwise use a number as the first argument or specify no optional argument.}}
        }{\str_if_eq:nnTF{#1}{e}
            {\IfValueTF{#2}
                {#2{\scalebox{0.65}[0.65]{$#3\lambda$}}}
                {\scalebox{0.65}[0.65]{$#3\lambda$}}
        }{\str_if_eq:nnTF{#1}{b}
                {\IfValueTF{#2}
                    {#2{#3\lambda}}
                    {{#3\lambda}}
         }{\IfValueTF{#2}
                    {#2{\scalebox{#1}[#1]{$#3\lambda$}}}
                    {\scalebox{#1}[#1]{$#3\lambda$}}
}}}}

\def\lg{
    \@ifgrave
        {
        \@ifstar
            {
            \@lg[e][\mathclap]
            }
            {
            \@lg[e]
            }
        }
        {    \@ifstar
             {
             \@lg[n][\mathclap]
             }
             {
             \@lg
             }
         }
     }

\newcommand{\xlg}{\@lg}

\NewDocumentCommand{\@mg}{O{0.85}oO{}}
{\str_if_eq:nnTF{#1}{n}
        {\IfValueTF{#2}
            {#2{\scalebox{0.85}[0.85]{$#3\mu$}}}
            {\errmessage{Too few arguments}\errhelp{If you specify `n' as the first argument, the second optional argument must be specified, otherwise use a number as the first argument or specify no optional argument.}}
        }{\str_if_eq:nnTF{#1}{e}
            {\IfValueTF{#2}
                {#2{\scalebox{0.65}[0.65]{$#3\mu$}}}
                {\scalebox{0.65}[0.65]{$#3\mu$}}
        }{\str_if_eq:nnTF{#1}{b}
                {\IfValueTF{#2}
                    {#2{#3\mu}}
                    {{#3\mu}}
         }{\str_if_eq:nnTF{#1}{t}
             {\IfValueTF{#2}
                 {#2{\scalebox{0.88}[0.88]{$#3\mu$}}}
                 {\scalebox{0.88}[0.88]{$#3\mu$}}
                 }{\IfValueTF{#2}
                    {#2{\scalebox{#1}[#1]{$#3\mu$}}}
                    {\scalebox{#1}[#1]{$#3\mu$}}
}}}}}

\newcommand{\mg}{
    \@ifgrave
        {
        \@ifstar
            {
            \@mg[e][\mathclap]
            }
            {
            \@mg[e]
            }
        }
        {    \@ifstar
             {
             \@mg[n][\mathclap]
             }
             {
             \@mg
             }
         }
     }

\newcommand{\xmg}{\@mg}
     
\NewDocumentCommand{\@ng}{O{0.85}oO{}}
{\str_if_eq:nnTF{#1}{n}
        {\IfValueTF{#2}
            {#2{\scalebox{0.85}[0.85]{$#3\nu$}}}
            {\errmessage{Too few arguments}\errhelp{If you specify `n' as the first argument, the second optional argument must be specified, otherwise use a number as the first argument or specify no optional argument.}}
        }{\str_if_eq:nnTF{#1}{e}
            {\IfValueTF{#2}
                {#2{\scalebox{0.65}[0.65]{$#3\nu$}}}
                {\scalebox{0.65}[0.65]{$#3\nu$}}
        }{\str_if_eq:nnTF{#1}{b}
                {\IfValueTF{#2}
                    {#2{#3\nu}}
                    {{#3\nu}}
         }{\IfValueTF{#2}
                    {#2{\scalebox{#1}[#1]{$#3\nu$}}}
                    {\scalebox{#1}[#1]{$#3\nu$}}
}}}}

\def\ng{
    \@ifgrave
        {
        \@ifstar
            {
            \@ng[e][\mathclap]
            }
            {
            \@ng[e]
            }
        }
        {    \@ifstar
             {
             \@ng[n][\mathclap]
             }
             {
             \@ng
             }
         }
     }

\newcommand{\xng}{\@ng}

\NewDocumentCommand{\@jg}{O{0.88}oO{}}
{\str_if_eq:nnTF{#1}{n}
        {\IfValueTF{#2}
            {#2{\scalebox{0.88}[0.88]{$#3\xi$}}}
            {\errmessage{Too few arguments}\errhelp{If you specify `n' as the first argument, the second optional argument must be specified, otherwise use a number as the first argument or specify no optional argument.}}
        }{\str_if_eq:nnTF{#1}{e}
            {\IfValueTF{#2}
                {#2{\scalebox{0.65}[0.65]{$#3\xi$}}}
                {\scalebox{0.65}[0.65]{$#3\xi$}}
        }{\str_if_eq:nnTF{#1}{b}
                {\IfValueTF{#2}
                    {#2{#3\xi}}
                    {{#3\xi}}
         }{\IfValueTF{#2}
                    {#2{\scalebox{#1}[#1]{$#3\xi$}}}
                    {\scalebox{#1}[#1]{$#3\xi$}}
}}}}

\newcommand{\jg}{
    \@ifgrave
        {
        \@ifstar
            {
            \@jg[e][\mathclap]
            }
            {
            \@jg[e]
            }
        }
        {    \@ifstar
             {
             \@jg[n][\mathclap]
             }
             {
             \@jg
             }
         }
     }

\newcommand{\xjg}{\@jg}

\NewDocumentCommand{\@og}{O{0.96}oO{}}
{\str_if_eq:nnTF{#1}{n}
        {\IfValueTF{#2}
            {#2{\scalebox{0.96}[0.96]{$#3\omicron$}}}
            {\errmessage{Too few arguments}\errhelp{If you specify `n' as the first argument, the second optional argument must be specified, otherwise use a number as the first argument or specify no optional argument.}}
        }{\str_if_eq:nnTF{#1}{e}
            {\IfValueTF{#2}
                {#2{\scalebox{0.72}[0.72]{$#3\omicron$}}}
                {\scalebox{0.72}[0.72]{$#3\omicron$}}
        }{\str_if_eq:nnTF{#1}{b}
                {\IfValueTF{#2}
                    {#2{#3\omicron}}
                    {{#3\omicron}}
         }{\IfValueTF{#2}
                    {#2{\scalebox{#1}[#1]{$#3\omicron$}}}
                    {\scalebox{#1}[#1]{$#3\omicron$}}
}}}}

\newcommand{\og}{
    \@ifgrave
        {
        \@ifstar
            {
            \@og[e][\mathclap]
            }
            {
            \@og[e]
            }
        }
        {    \@ifstar
             {
             \@og[n][\mathclap]
             }
             {
             \@og
             }
         }
     }

\newcommand{\xog}{\@og}

\NewDocumentCommand{\@pg}{O{0.85}oO{}}
{\str_if_eq:nnTF{#1}{n}
        {\IfValueTF{#2}
            {#2{\scalebox{0.85}[0.85]{$#3\pi$}}}
            {\errmessage{Too few arguments}\errhelp{If you specify `n' as the first argument, the second optional argument must be specified, otherwise use a number as the first argument or specify no optional argument.}}
        }{\str_if_eq:nnTF{#1}{e}
            {\IfValueTF{#2}
                {#2{\scalebox{0.65}[0.65]{$#3\pi$}}}
                {\scalebox{0.65}[0.65]{$#3\pi$}}
        }{\str_if_eq:nnTF{#1}{b}
                {\IfValueTF{#2}
                    {#2{#3\pi}}
                    {{#3\pi}}
         }{\IfValueTF{#2}
                    {#2{\scalebox{#1}[#1]{$#3\pi$}}}
                    {\scalebox{#1}[#1]{$#3\pi$}}
}}}}

\newcommand{\pg}{
    \@ifgrave
        {
        \@ifstar
            {
            \@pg[e][\mathclap]
            }
            {
            \@pg[e]
            }
        }
        {    \@ifstar
             {
             \@pg[n][\mathclap]
             }
             {
             \@pg
             }
         }
     }
     
\newcommand{\xpg}{\@pg}

\NewDocumentCommand{\@rg}{O{0.85}oO{}}
{\str_if_eq:nnTF{#1}{n}
        {\IfValueTF{#2}
            {#2{\scalebox{0.85}[0.85]{$#3\rho$}}}
            {\errmessage{Too few arguments}\errhelp{If you specify `n' as the first argument, the second optional argument must be specified, otherwise use a number as the first argument or specify no optional argument.}}
        }{\str_if_eq:nnTF{#1}{e}
            {\IfValueTF{#2}
                {#2{\scalebox{0.65}[0.65]{$#3\rho$}}}
                {\scalebox{0.65}[0.65]{$#3\rho$}}
        }{\str_if_eq:nnTF{#1}{b}
                {\IfValueTF{#2}
                    {#2{#3\rho}}
                    {{#3\rho}}
         }{\str_if_eq:nnTF{#1}{t}
             {\IfValueTF{#2}
                 {#2{\scalebox{0.88}[0.88]{$#3\rho$}}}
                 {\scalebox{0.88}[0.88]{$#3\rho$}}
                 }{\IfValueTF{#2}
                    {#2{\scalebox{#1}[#1]{$#3\rho$}}}
                    {\scalebox{#1}[#1]{$#3\rho$}}
}}}}}

\newcommand{\rg}{
    \@ifgrave
        {
        \@ifstar
            {
            \@rg[e][\mathclap]
            }
            {
            \@rg[e]
            }
        }
        {    \@ifstar
             {
             \@rg[n][\mathclap]
             }
             {
             \@rg
             }
         }
     }
     
\newcommand{\xrg}{\@rg}

\NewDocumentCommand{\@sg}{O{0.85}oO{}}
{\str_if_eq:nnTF{#1}{n}
        {\IfValueTF{#2}
            {#2{\scalebox{0.85}[0.85]{$#3\sigma$}}}
            {\errmessage{Too few arguments}\errhelp{If you specify `n' as the first argument, the second optional argument must be specified, otherwise use a number as the first argument or specify no optional argument.}}
        }{\str_if_eq:nnTF{#1}{e}
            {\IfValueTF{#2}
                {#2{\scalebox{0.65}[0.65]{$#3\sigma$}}}
                {\scalebox{0.65}[0.65]{$#3\sigma$}}
        }{\str_if_eq:nnTF{#1}{b}
                {\IfValueTF{#2}
                    {#2{#3\sigma}}
                    {{#3\sigma}}
         }{\IfValueTF{#2}
                    {#2{\scalebox{#1}[#1]{$#3\sigma$}}}
                    {\scalebox{#1}[#1]{$#3\sigma$}}
}}}}

\newcommand{\sg}{
    \@ifgrave
        {
        \@ifstar
            {
            \@sg[e][\mathclap]
            }
            {
            \@sg[e]
            }
        }
        {    \@ifstar
             {
             \@sg[n][\mathclap]
             }
             {
             \@sg
             }
         }
     }
     
\newcommand{\xsg}{\@sg}

\NewDocumentCommand{\@tg}{O{0.85}oO{}}
{\str_if_eq:nnTF{#1}{n}
        {\IfValueTF{#2}
            {#2{\scalebox{0.85}[0.85]{$#3\tau$}}}
            {\errmessage{Too few arguments}\errhelp{If you specify `n' as the first argument, the second optional argument must be specified, otherwise use a number as the first argument or specify no optional argument.}}
        }{\str_if_eq:nnTF{#1}{e}
            {\IfValueTF{#2}
                {#2{\scalebox{0.65}[0.65]{$#3\tau$}}}
                {\scalebox{0.65}[0.65]{$#3\tau$}}
        }{\str_if_eq:nnTF{#1}{b}
                {\IfValueTF{#2}
                    {#2{#3\tau}}
                    {{#3\tau}}
         }{\IfValueTF{#2}
                    {#2{\scalebox{#1}[#1]{$#3\tau$}}}
                    {\scalebox{#1}[#1]{$#3\tau$}}
}}}}

\newcommand{\tg}{
    \@ifgrave
        {
        \@ifstar
            {
            \@tg[e][\mathclap]
            }
            {
            \@tg[e]
            }
        }
        {    \@ifstar
             {
             \@tg[n][\mathclap]
             }
             {
             \@tg
             }
         }
     }
     
\newcommand{\xtg}{\@tg}

\NewDocumentCommand{\@ug}{O{0.85}oO{}}
{\str_if_eq:nnTF{#1}{n}
        {\IfValueTF{#2}
            {#2{\scalebox{0.85}[0.85]{$#3\upsilon$}}}
            {\errmessage{Too few arguments}\errhelp{If you specify `n' as the first argument, the second optional argument must be specified, otherwise use a number as the first argument or specify no optional argument.}}
        }{\str_if_eq:nnTF{#1}{e}
            {\IfValueTF{#2}
                {#2{\scalebox{0.65}[0.65]{$#3\upsilon$}}}
                {\scalebox{0.65}[0.65]{$#3\upsilon$}}
        }{\str_if_eq:nnTF{#1}{b}
                {\IfValueTF{#2}
                    {#2{#3\upsilon}}
                    {{#3\upsilon}}
         }{\IfValueTF{#2}
                    {#2{\scalebox{#1}[#1]{$#3\upsilon$}}}
                    {\scalebox{#1}[#1]{$#3\upsilon$}}
}}}}

\newcommand{\ug}{
    \@ifgrave
        {
        \@ifstar
            {
            \@ug[e][\mathclap]
            }
            {
            \@ug[e]
            }
        }
        {    \@ifstar
             {
             \@ug[n][\mathclap]
             }
             {
             \@ug
             }
         }
     }
     
\newcommand{\xug}{\@ug}

\NewDocumentCommand{\@fg}{O{0.85}oO{}}
{\str_if_eq:nnTF{#1}{n}
        {\IfValueTF{#2}
            {#2{\scalebox{0.85}[0.85]{$#3\varphi$}}}
            {\errmessage{Too few arguments}\errhelp{If you specify `n' as the first argument, the second optional argument must be specified, otherwise use a number as the first argument or specify no optional argument.}}
        }{\str_if_eq:nnTF{#1}{e}
            {\IfValueTF{#2}
                {#2{\scalebox{0.65}[0.65]{$#3\varphi$}}}
                {\scalebox{0.65}[0.65]{$#3\varphi$}}
        }{\str_if_eq:nnTF{#1}{b}
                {\IfValueTF{#2}
                    {#2{#3\varphi}}
                    {{#3\varphi}}
         }{\str_if_eq:nnTF{#1}{t}
             {\IfValueTF{#2}
                 {#2{\scalebox{0.88}[0.88]{$#3\varphi$}}}
                 {\scalebox{0.88}[0.88]{$#3\varphi$}}
                 }{\IfValueTF{#2}
                    {#2{\scalebox{#1}[#1]{$#3\varphi$}}}
                    {\scalebox{#1}[#1]{$#3\varphi$}}
}}}}}

\newcommand{\fg}{
    \@ifgrave
        {
        \@ifstar
            {
            \@fg[e][\mathclap]
            }
            {
            \@fg[e]
            }
        }
        {    \@ifstar
             {
             \@fg[n][\mathclap]
             }
             {
             \@fg
             }
         }
     }
     
\newcommand{\xfg}{\@fg}

\NewDocumentCommand{\@xg}{O{0.85}oO{}}
{\str_if_eq:nnTF{#1}{n}
        {\IfValueTF{#2}
            {#2{\scalebox{0.85}[0.85]{$#3\chi$}}}
            {\errmessage{Too few arguments}\errhelp{If you specify `n' as the first argument, the second optional argument must be specified, otherwise use a number as the first argument or specify no optional argument.}}
        }{\str_if_eq:nnTF{#1}{e}
            {\IfValueTF{#2}
                {#2{\scalebox{0.65}[0.65]{$#3\chi$}}}
                {\scalebox{0.65}[0.65]{$#3\chi$}}
        }{\str_if_eq:nnTF{#1}{b}
                {\IfValueTF{#2}
                    {#2{#3\chi}}
                    {{#3\chi}}
         }{\str_if_eq:nnTF{#1}{t}
             {\IfValueTF{#2}
                 {#2{\scalebox{0.88}[0.88]{$#3\chi$}}}
                 {\scalebox{0.88}[0.88]{$#3\chi$}}
                 }{\IfValueTF{#2}
                    {#2{\scalebox{#1}[#1]{$#3\chi$}}}
                    {\scalebox{#1}[#1]{$#3\chi$}}
}}}}}

\newcommand{\xg}{
    \@ifgrave
        {
        \@ifstar
            {
            \@xg[e][\mathclap]
            }
            {
            \@xg[e]
            }
        }
        {    \@ifstar
             {
             \@xg[n][\mathclap]
             }
             {
             \@xg
             }
         }
     }
     
\newcommand{\xxg}{\@xg}

\NewDocumentCommand{\@yg}{O{0.85}oO{}}
{\str_if_eq:nnTF{#1}{n}
        {\IfValueTF{#2}
            {#2{\scalebox{0.85}[0.85]{$#3\psi$}}}
            {\errmessage{Too few arguments}\errhelp{If you specify `n' as the first argument, the second optional argument must be specified, otherwise use a number as the first argument or specify no optional argument.}}
        }{\str_if_eq:nnTF{#1}{e}
            {\IfValueTF{#2}
                {#2{\scalebox{0.65}[0.65]{$#3\psi$}}}
                {\scalebox{0.65}[0.65]{$#3\psi$}}
        }{\str_if_eq:nnTF{#1}{b}
                {\IfValueTF{#2}
                    {#2{#3\psi}}
                    {{#3\psi}}
         }{\str_if_eq:nnTF{#1}{t}
             {\IfValueTF{#2}
                 {#2{\scalebox{0.88}[0.88]{$#3\psi$}}}
                 {\scalebox{0.88}[0.88]{$#3\psi$}}
                 }{\IfValueTF{#2}
                    {#2{\scalebox{#1}[#1]{$#3\psi$}}}
                    {\scalebox{#1}[#1]{$#3\psi$}}
}}}}}

\newcommand{\yg}{
    \@ifgrave
        {
        \@ifstar
            {
            \@yg[e][\mathclap]
            }
            {
            \@yg[e]
            }
        }
        {    \@ifstar
             {
             \@yg[n][\mathclap]
             }
             {
             \@yg
             }
         }
     }
     
\newcommand{\xyg}{\@yg}

\NewDocumentCommand{\@wg}{O{0.85}oO{}}
{\str_if_eq:nnTF{#1}{n}
        {\IfValueTF{#2}
            {#2{\scalebox{0.85}[0.85]{$#3\omega$}}}
            {\errmessage{Too few arguments}\errhelp{If you specify `n' as the first argument, the second optional argument must be specified, otherwise use a number as the first argument or specify no optional argument.}}
        }{\str_if_eq:nnTF{#1}{e}
            {\IfValueTF{#2}
                {#2{\scalebox{0.65}[0.65]{$#3\omega$}}}
                {\scalebox{0.65}[0.65]{$#3\omega$}}
        }{\str_if_eq:nnTF{#1}{b}
                {\IfValueTF{#2}
                    {#2{#3\omega}}
                    {{#3\omega}}
         }{\IfValueTF{#2}
                    {#2{\scalebox{#1}[#1]{$#3\omega$}}}
                    {\scalebox{#1}[#1]{$#3\omega$}}
}}}}

\newcommand{\wg}{
    \@ifgrave
        {
        \@ifstar
            {
            \@wg[e][\mathclap]
            }
            {
            \@wg[e]
            }
        }
        {    \@ifstar
             {
             \@wg[n][\mathclap]
             }
             {
             \@wg
             }
         }
     }
     
\newcommand{\xwg}{\@wg}
\def\@ifgrave#1{\@ifnextchar{`}{\@firstoftwo{#1}}}
\def\@ifstar#1{\kernel@ifnextchar *{\@firstoftwo{#1}}}
\newlength\mw@frraise
\newlength\mw@frlower
\NewDocumentCommand{\ifr}{mmO{1.5}O{-1.5}}
{\frac{\raisebox{#3pt}{$\scriptstyle #1$}}{\raisebox{#4pt}{$\scriptstyle #2$}}}
\NewDocumentCommand{\@fr}{mmO{-1pt}O{0pt}}
{\dfrac{\raisebox{#3}{$\displaystyle #1$}}{\raisebox{#4}{$\displaystyle #2$}}}
\NewDocumentCommand{\@dbfr}{mmO{-1pt}O{0pt}}
{\dbfr{\raisebox{#3}{$#1$}}{\raisebox{#4}{$#2$}}}
\NewDocumentCommand{\dcfr}{mmO{-1pt}O{0pt}D**{0.4}}
{\scalebox{#5}[#5]{$\fr{\raisebox{#3}{$#1$}}{\raisebox{#4}{$#2$}}$}}

\newcommand{\dbfr}[2]{\dsib{\fr{#1}{#2}}}

\newcommand{\dsib}[1]{\scb{0.5}[0.5]{$#1$}}
\newcommand{\scb}{\scalebox}
\NewDocumentCommand{\xfr}{mmood**}
{\begingroup
\IfValueTF{#3}{
    \str_if_eq:nnTF{#3}{def}
        {\setlength{\mw@frraise}{1.5pt}}
        {\setlength{\mw@frraise}{#3pt}}
    }{}
\IfValueTF{#4}{\setlength{\mw@frlower}{#4pt}}{}
\IfValueTF{#5}{\def\mw@xfr@scale{#5}}{}
\mathchoice%
{\advance\mw@frraise by -2.5pt
\advance\mw@frlower by 1.5pt
\@fr{#1}{#2}[\mw@frraise][\mw@frlower]}%
{\ifr{#1}{#2}[\strip@pt\mw@frraise][\strip@pt\mw@frlower]}%
{
\advance\mw@frraise by -2.5pt
\advance\mw@frlower by 1.5pt
\@dbfr{#1}{#2}[\mw@frraise][\mw@frlower]
}%
{
\advance\mw@frraise by -2.5pt
\advance\mw@frlower by 1.5pt
\dcfr{#1}{#2}[\mw@frraise][\mw@frlower]*\mw@xfr@scale*
}%
\endgroup}
\NewDocumentCommand{\setxfr}{od()d**}
{\IfValueTF{#1}{\setlength{\mw@frraise}{#1pt}}{}
\IfValueTF{#2}{\setlength{\mw@frlower}{#2pt}}{}
\IfValueTF{#3}{\def\mw@xfr@scale{#3}}
}
\setxfr[1.5](-1.5)*0.4*\relax
\newcommand{\delim}[3][lr]{
{
\str_case:nnF { #1 }
    {
    { b }{\def\lsize{\big}\def\rsize{\big}}
    { B }{\def\lsize{\Big}\def\rsize{\Big}}
    { x }{\def\lsize{\bigg}\def\rsize{\bigg}}
    { X }{\def\lsize{\Bigg}\def\rsize{\Bigg}}
    { lr }{\def\lsize{\left}\def\rsize{\right}    \lrtrue}
    { bB }{\def\lsize{\big}\def\rsize{\Big}}
    { Bb }{\def\lsize{\Big}\def\rsize{\big}}
    { xb }{\def\lsize{\bigg}\def\rsize{\big}}
    { xB }{\def\lsize{\bigg}\def\rsize{\Big}}
    { bx }{\def\lsize{\big}\def\rsize{\bigg}}
    { Bx }{\def\lsize{\Big}\def\rsize{\bigg}}
    { bX }{\def\lsize{\big}\def\rsize{\Bigg}}
    { BX }{\def\lsize{\Big}\def\rsize{\Bigg}}
    { xX }{\def\lsize{\bigg}\def\rsize{\Bigg}}
    { Xb }{\def\lsize{\Bigg}\def\rsize{\big}}
    { XB }{\def\lsize{\Bigg}\def\rsize{\Big}}
    { Xx }{\def\lsize{\Bigg}\def\rsize{\bigg}}
    {n}{\def\lsize{}\def\rsize{}}
    {bn}{\def\lsize{\big}\def\rsize{}}
    {Bn}{\def\lsize{\Big}\def\rsize{}}
    {xn}{\def\lsize{\bigg}\def\rsize{}}
    {Xn}{\def\lsize{\Bigg}\def\rsize{}}
    {nb}{\def\lsize{}\def\rsize{\big}}
    {nB}{\def\lsize{}\def\rsize{\Big}}
    {nx}{\def\lsize{}\def\rsize{\bigg}}
    {nX}{\def\lsize{}\def\rsize{\Bigg}}
}
{\msg_fatal:nnn{mworksx}{invalid}{delimiter~size}}
\str_case:nnF { #2 }
    {
    { s }{\def\ldel{[}\def\rdel{]}}
    { r }{\def\ldel{(}\def\rdel{)}}
    { b }{\def\ldel{\lbrace}\def\rdel{\rbrace}}
    { v }{\def\ldel{|}\def\rdel{|}}
    { a }{\def\ldel{\langle}\def\rdel{\rangle}}
    { dv }{\def\ldel{\|}\def\rdel{\|}}
    { rs }{\def\ldel{(}\def\rdel{]}}
    { sr }{\def\ldel{[}\def\rdel{)}}
    { rb }{\def\ldel{(}\def\rdel{\rbrace}}
    { sb }{\def\ldel{[}\def\rdel{\rbrace}}
    { br }{\def\ldel{\lbrace}\def\rdel{)}}
    { bs }{\def\ldel{\lbrace}\def\rdel{]}}
    { ra }{\def\ldel{(}\def\rdel{\rangle}}
    { sa }{\def\ldel{[}\def\rdel{\rangle}}
    { ba }{\def\ldel{\lbrace}\def\rdel{\rangle}}
    { ar }{\def\ldel{\langle}\def\rdel{)}}
    { as }{\def\ldel{\langle}\def\rdel{]}}
    { ab }{\def\ldel{\langle}\def\rdel{\rbrace}}
    { rv }{\def\ldel{(}\def\rdel{|}}
    { sv }{\def\ldel{[}\def\rdel{|}}
    { bv }{\def\ldel{\lbrace}\def\rdel{|}}
    { av }{\def\ldel{\langle}\def\rdel{|}}
    { vr }{\def\ldel{|}\def\rdel{)}}
    { vs }{\def\ldel{|}\def\rdel{]}}
    { vb }{\def\ldel{|}\def\rdel{\rbrace}}
    { va }{\def\ldel{|}\def\rdel{\rangle}}
    { rdv }{\def\ldel{(}\def\rdel{\|}}
    { sdv }{\def\ldel{[}\def\rdel{\|}}
    { bdv }{\def\ldel{\lbrace}\def\rdel{\|}}
    { adv }{\def\ldel{\langle}\def\rdel{\|}}
    { vdv }{\def\ldel{|}\def\rdel{\|}}
    { dvr }{\def\ldel{\|}\def\rdel{)}}
    { dvs }{\def\ldel{\|}\def\rdel{]}}
    { dvb }{\def\ldel{\|}\def\rdel{\rbrace}}
    { dva }{\def\ldel{\|}\def\rdel{\rangle}}
    { dvv }{\def\ldel{\|}\def\rdel{|}}
    { lr }{\def\ldel{(}\def\rdel{.}}
    { rr }{\def\ldel{)}\def\rdel{.}}
    { rir }{\def\ldel{.}\def\rdel{)}}
    { ls }{\def\ldel{[}\def\rdel{.}}
    { rs }{\def\ldel{]}\def\rdel{.}}
    { ris }{\def\ldel{.}\def\rdel{]}}
    { lb }{\def\ldel{\lbrace}\def\rdel{.}}
    { rb }{\def\ldel{\rbrace}\def\rdel{.}}
    { rib }{\def\ldel{.}\def\rdel{\rbrace}}
    { la }{\def\ldel{\langle}\def\rdel{.}}
    { ra }{\def\ldel{\rangle}\def\rdel{.}}
    { ria }{\def\ldel{.}\def\rdel{\rangle}}
    { ov }{\def\ldel{|}\def\rdel{.}}
    { odv }{\def\ldel{\|}\def\rdel{.}}
    { ssi }{\invd@limstrue
        \iflr
            \def\ldel{[}\def\rdel{[}\invd@limsfalse
        \else
            \def\ldel{\lsize[}\def\rdel{\mathclose{\rsize[}}
        \fi}
    { sis }{\invd@limstrue
        \iflr
            \def\ldel{]}\def\rdel{]}\invd@limsfalse
        \else
            \def\ldel{\mathopen{\lsize]}}\def\rdel{\rsize]}
        \fi}
    { sisi }{\invd@limstrue
        \iflr
            \def\ldel{]}\def\rdel{[}\invd@limsfalse
        \else
            \def\ldel{\mathopen{\lsize]}}\def\rdel{\mathclose{\rsize[}}
        \fi}
}
{\msg_fatal:nnn{mworksx}{invalid}{delimiter}}
\ifx\rdel.
    \lsize\ldel
\else
    \ifinvd@lims
        \ldel #3 \rdel
    \else
        \lsize\ldel #3 \rsize\rdel
    \fi
\fi
}
}
\NewDocumentCommand{\ux}{sO{}mm}
{\IfBooleanTF{#1}
    {\overset{\xstack*[#2]{#3}}{#4}}
    {\underset{\xstack[#2]{#3}}{#4}}
}
\NewDocumentCommand{\ua}{sO{}mm}
{\IfBooleanTF{#1}
    {\overset{\arstack*[#2]{#3}}{#4}}
    {\underset{\arstack[#2]{#3}}{#4}}
}
\NewDocumentCommand{\pux}{sO{}D**{8em}mm}
{\IfBooleanTF{#1}
    {\overset{\xstack*[#2]{\parbox{#3}{\centering\linespread{1.0}\selectfont\scriptsize #4}}}{#5}}
    {\underset{\xstack[#2]{\parbox{#3}{\centering\linespread{1.0}\selectfont\scriptsize #4}}}{#5}}
}
\NewDocumentCommand{\pua}{sO{}D**{8em}mm}
{\IfBooleanTF{#1}
    {\overset{\arstack*[#2]{\parbox{#3}{\centering\linespread{1.0}\selectfont\scriptsize #4}}}{#5}}
    {\underset{\arstack[#2]{\parbox{#3}{\centering\linespread{1.0}\selectfont\scriptsize #4}}}{#5}}
}
\NewDocumentCommand{\puxtrick}{sO{}D**{8em}m}
{\IfBooleanTF{#1}
     {\xstack*[#2]{\parbox{#3}{\centering\linespread{1.0}\selectfont\scriptsize #4}}}
     {\xstack[#2]{\parbox{#3}{\centering\linespread{1.0}\selectfont\scriptsize #4}}}
}
\NewDocumentCommand{\puatrick}{sO{}D**{8em}m}
{\IfBooleanTF{#1}
     {\arstack*[#2]{\parbox{#3}{\centering\linespread{1.0}\selectfont\scriptsize #4}}}
     {\arstack[#2]{\parbox{#3}{\centering\linespread{1.0}\selectfont\scriptsize #4}}}
}
\NewDocumentCommand{\case}{mO{@{\,}cl@{}}}
{\left\lbrace\mat{#2}#1\emat\right.}

\newenvironment{sistema}[1][@{}l@{}]%
{\left \lbrace \begin{array}{#1}}%
{\end{array}\right.}
\RequirePackage{letltxmacro}
\LetLtxMacro{\oldsqrt}{\sqrt} 
\renewcommand{\sqrt}[1][\ ]{%
  \def\DHLindex{#1}\mathpalette\DHLhksqrt}
\def\DHLhksqrt#1#2{%
  \setbox0=\hbox{$#1\expandafter\oldsqrt\expandafter[\DHLindex]{#2\,}$}\dimen0=\ht0
  \advance\dimen0-0.2\ht0
  \setbox2=\hbox{\vrule height\ht0 depth -\dimen0}%
  {\box0\lower0.71pt\box2}}
\NewDocumentCommand{\kcref}{mo}{
    \expandafter\ifx\csname r@#1\endcsname\relax
    \protect\G@refundefinedtrue
   \textbf{(??)}
   \msg_warning:nnn{mworksx}{ref-undef}{#1}
   \else
   \IfValueTF{#2}
       {\def\@kind{#2}}
       {\edef\@labelfour{\expandafter\expandafter\expandafter\@fourthoffour\csname r@#1\endcsname}
   \edef\@argfour{\expandafter\@dotcar\@labelfour\relax.\relax\@nil}
   \@ifundefined{mw@\@argfour name}
       {\msg_warning:nnn{mworksx}{thm-kind-undef}{\@argfour}
       \def\@kind{\hspace{-\mw@sp}}}
       {\edef\@kind{\csname mw@\@argfour name\endcsname}}}
   \hyperref[#1]{\refstyle{\@kind\ \ref*{#1}}}
   \edef\@defaulttype{\csname default\@argfour type\endcsname}
   \edef\@fixtype{\noexpand\setcnttype{\@argfour}{\@defaulttype}}
   \@fixtype
   \fi
}
\NewDocumentCommand{\kprepnref}{mo}{
    \expandafter\ifx\csname r@#1\endcsname\relax
    \protect\G@refundefinedtrue
   \textbf{(??)}
   \msg_warning:nnn{mworksx}{ref-undef}{#1}
   \else
   \IfValueTF{#2}
       {\def\@kind{#2}}
       {\edef\@labelfour{\expandafter\expandafter\expandafter\@fourthoffour\csname r@#1\endcsname}
   \edef\@argfour{\expandafter\@dotcar\@labelfour\relax.\relax\@nil}
   \@ifundefined{mw@\@argfour name}
       {\msg_warning:nnn{mworksx}{thm-kind-undef}{\@argfour}
       \def\@kind{\hspace{-\mw@sp}}}
       {\edef\@kind{\csname mw@\@argfour name\endcsname}}}
   \hyperref[#1]{\refstyle{\@kind\ \refprep\ \nameref*{#1}}}
   \fi
}
\NewDocumentCommand{\kcolonnref}{mo}{
    \expandafter\ifx\csname r@#1\endcsname\relax
    \protect\G@refundefinedtrue
   \textbf{(??)}
   \msg_warning:nnn{mworksx}{ref-undef}{#1}
   \else
   \IfValueTF{#2}
       {\def\@kind{#2}}
       {\edef\@labelfour{\expandafter\expandafter\expandafter\@fourthoffour\csname r@#1\endcsname}
   \edef\@argfour{\expandafter\@dotcar\@labelfour\relax.\relax\@nil}
   \@ifundefined{mw@\@argfour name}
       {\msg_warning:nnn{mworksx}{thm-kind-undef}{\@argfour}
       \def\@kind{\hspace{-\mw@sp}}}
       {\edef\@kind{\csname mw@\@argfour name\endcsname}}}
   \hyperref[#1]{\refstyle{\@kind\refpunct\ \nameref*{#1}}}
   \fi
}
\msg_new:nnnn{mworksx}{too-few-args}{Too ~ few ~ arguments.}{The ~ macro ~ \string#1 ~ requires ~ at ~ least ~ one ~ of ~ its ~ arguments ~ to ~ have ~ a ~ value. ~ Please ~ supply ~ it.}
\msg_new:nnn{mworksx}{save-no-inpenc}{You\ are\ using\ inputenc\ with\ the\ save\ option\ from\ this\ package.\ Be\ careful\ with\ characters\ with\ tilde\ or\ double\ acute\ accent\ above,\ as\ this\ package\ may\ render\ them\ impossible\ to\ typeset,\ or\ else\ fix\ the\ issue\ by\ the\ inpenc\ option\ of\ this\ package.\ Note\ also\ that\ the\ command\ \string\ss\ is\ redefined\ to\ be\ the\ spadesuit\ but\ to\ switch\ back\ to\ the\ German\ Scharfeses\ (ß)\ if\ starred.\ If\ there\ is\ any\ other\ trouble,\ which\ vanishes\ if\ you\ remove\ the\ loading\ of\ mworksx,\ contact\ me\ at\ michelegor@hotmail.com.\ Also,\ the\ utf8x\ option\ doesn't\ interact\ well\ with\ mworksx,\ so\ I\ suggest\ you\ do\ NOT\ use\ it\ with\ this\ package.\ Thanks\ for\ your\ attention\ and\ I\ hope\ this\ package\ serves\ you\ well\ :).}
\msg_new:nnnn{mworksx}{invalid}{What~exactly~did~you~mean?}{
I~don't~get~what~delimiter~size~you~mean~with~this.~Please~review~your~#1~specifications.}
\msg_new:nnnn{mworksx}{not-read}{You\ haven't\ read\ the\ documentation\ yet!}{This\ used\ to\ be\ an\ error.\ With\ the\ new\ version\ of\ the\ package,\ new\ warnings\ make\ it,\ in\ my\ opinion,\ not\ strictly\ necessary\ to\ read\ it.\ However,\ I\ strongly\ recommend\ that\ you\ read\ it,\ both\ to\ familiarize\ with\ the\ commands\ and\ to\ know\ the\ reason\ for\ this\ warning.\ Thanks\ for\ using\ mworksx.}
\msg_new:nnnn{mworksx}{unknown-lang}{Unknown~language~option~for~something~in~this~package.}{Couldn't~interpret~given~language~options.~#1~left~empty.}
\msg_new:nnnn{mworksx}{invalid-ctr}{Illegal\ argument\ in\ counter\ definition.}{The\ third\ argument\ must\ be\ either\ c,\ s\ or\ x,\ for\ chapter,\ section\ and\ subsection\ respectively,\ or\ n\ for\ none\ of\ those\ which\ means\ just\ the\ counter\ in\ parameter\ 2,\ or\ a\ to\ let\ me\ choose\ the\ counters\ automatically,\ or\ @@\ for\ no\ counters\ at\ all.}
\msg_new:nnn{mworksx}{label-undefined}{The ~ label ~ #2 ~ is ~ not ~ defined ~ or ~ cannot ~ be ~ found.}
\msg_new:nnn{mworksx}{ref-undef}{Reference~#1~undefined~on~page~\thepage.}
\msg_new:nnn{mworksx}{cnt-type-undef}{Counter~type~#1~undefined.~The~known~types~are~c,~s,~x,~a,~n,~nc~and~@@.~Please~choose~one~of~these.}
\msg_new:nnn{mworksx}{thm-kind-undef}{Theorem~kind~#1~not~found.~Please~define~\string\mw@ #1name.~Printing~nothing.}
\msg_new:nnnn{mworksx}{unknown-opt-value}{Unknown~value~given~for~option~#1.}{I~can't~understand~what~you~wanted~when~you~gave~the~value~#2~for~option~#1.~Please~review~and~correct~your~choice,~or~ignore~the~warning~and~stick~to~the~default~#3.}
\msg_new:nnnn{mworksx}{no-dfint}{Double-dash~integral~unavailable}{Fixing~option~fint~at~slanted~doesn't~give~integral~symbol~with~double~strikethrough.~If~that~is~needed,~set~fint~to~straight~or~fintmint~to~true.}
\NewDocumentCommand{\eval}{O{B}mo}
{\IfValueTF{#3}
    {\sthat[#1]_{#2}^{#3}}
    {\sthat[#1]_{#2}}
}
\newcommand{\mw@ptmx}{all=upnewtx}
\input newtxstuff.tex
\ExplSyntaxOff
\newcommand{\xint}{\int\limits}
\NewDocumentCommand{\xints}{mmod**}
{\IfValueTF{#3}{\mw@intlow@space=#3}{}
\IfValueTF{#4}{\mw@inthigh@space=#4}{}
\xint_{#1\mskip\mw@intlow@space}^{\mskip\mw@inthigh@space #2}}
\newmuskip\mw@intlow@space
\newmuskip\mw@inthigh@space
\newcommand{\sps}{\@ifstar{\supset}{\!\supset\!}}
\newcommand{\spse}{\@ifstar{\supseteq}{\!\supseteq\!}}
\newcommand{\sbs}{\@ifstar{\subset}{\!\subset\!}}
\newcommand{\sbse}{\@ifstar{\subseteq}{\!\subseteq\!}}
\newcommand{\nsps}{\@ifstar{\centernot\supset}{\!\centernot\supset\!}}
\newcommand{\nspse}{\@ifstar{\centernot\supseteq}{\!\centernot\supseteq\!}}
\newcommand{\nsbs}{\@ifstar{\centernot\subset}{\!\centernot\subset\!}}
\newcommand{\nsbse}{\@ifstar{\centernot\subseteq}{\!\centernot\subseteq\!}}
\let\basin\in
\let\basni\ni
\newcommand{\nin}{\@ifstar{\centernot\basin}{\!\centernot\basin\!}}
\newcommand{\nni}{\@ifstar{\centernot\basni}{\!\centernot\basni\!}}
\renewcommand{\in}{\@ifstar{\basin}{\!\basin\!}}
\renewcommand{\ni}{\@ifstar{\basni}{\!\basni\!}}

    \def\L{\mathbb L}
    
\newcommand{\N}{\mathbb N}
    
    \def\P{\mathbb P}

\newcommand{\R}{\mathbb R}
    
    \def\T{\mathbb T}

\newcommand{\two}{\@ifstar{\hsp{2pt}}{\hsp{-2pt}}}
\newcommand{\gettwo}{\let\2\two}
\gettwo
\newcommand{\xstack}{\@ifstar{\xst@ck}{\@xstack}}
\newcommand{\xst@ck}[2][]{\mcl{\sstack{#2\\#1|}}}
\newcommand{\@xstack}[2][]{\mcl{\sstack{#1|\\#2}}}

\newcommand{\hra}{\hookrightarrow}

\let\basto\to
\renewcommand{\to}{\@ifstar{\!\basto\!}{\basto}}

\newcommand{\VA}{\forall}
\newcommand{\8}{\infty}
\newcommand{\0}{\@ifstar{\emptyset}{\varnothing}}
\newcommand{\mat}{\begin{array}}
\newcommand{\emat}{\end{array}}
\newcommand{\x}{\times}
\newcommand{\ssm}{\smallsetminus}
\newcommand{\sthat}[1][B]{\delim[#1]{ov}{}}

\def\re@DeclareMathSymbol#1#2#3#4{%
    \let#1=\relax
    \DeclareMathSymbol{#1}{#2}{#3}{#4}}
\DeclareFontEncoding{LMX}{}{}
  \DeclareFontSubstitution{LMX}{ntxexx}{m}{n}
  \DeclareFontFamily{LMX}{ntxexx}{}
  \DeclareFontShape{LMX}{ntxexx}{m}{n}{<->ntxexx}{}
  \DeclareFontShape{LMX}{ntxexx}{b}{n}{<->ntxbexx}{}
  \DeclareFontShape{LMX}{ntxexx}{bx}{n}{<->ssub ntxexx/b/n}{}
  \DeclareSymbolFont{largesymbolsA}{LMX}{ntxexx}{m}{n}
  \SetSymbolFont{largesymbolsA}{bold}{LMX}{ntxexx}{b}{n}
  \re@DeclareMathSymbol{\intop}{\mathop}{largesymbolsA}{"52}
\newcommand{\@fourthoffour}[4]{#4}
\def\@dotcar#1.#2\@nil{#1}
\def\@dotcbr#1.#2\@nil{#2}
\makeatother

\newcommand{\Dg}{\Delta}

\newcommand{\Fg}{\Phi}

\newcommand{\per}{\cdot}
\newcommand{\rad}{\sqrt}
\newcommand{\ang}[1][lr]{\delim[#1]{a}}
\newcommand{\pa}[1][lr]{\delim[#1]{r}}
\newcommand{\sq}[1][lr]{\delim[#1]{s}}

\newcommand{\abs}[1][lr]{\delim[#1]{v}}
\newcommand{\norm}[1][lr]{\delim[#1]{dv}}

\newcommand{\mcl}{\mathclap}
\newcommand{\sstack}{\substack}
\newcommand{\ubar}{\underline}
\newcommand{\lbar}{\overline}
\newcommand{\pint}[1]{\mathaccent23{#1}}

\newcommand{\s}{\mathcal}

\newcommand{\Cinf}{\s{C}^\8}
\newcommand{\refstyle}[1]{{\bfseries\color{blue}#1}}

\newcommand{\lsim}{\lesssim}
\midv
\renewcommand{\thesection}{\arabic{section}}

\makeatletter
\newcommand{\frl}{\@ifstar{\@frl}{\frl@}}
\newcommand{\@frl}[1]{(-\Dg)^{\xfr{#1}{2}}}
\newcommand{\frl@}[1]{(-\Dg)^{#1}}
\makeatother
\definethm{probl}{Problem}[s]
\definethm{psprobl}{Pseudo-problem}[s]
\definethm{typo}{Probable typo/misprint}[s]
\definethm{genoss}{Navier-Stokes vs. Euler Remark}[s]
\definethm{conj}{Conjecture}[s]

\makeatletter
\newcommand{\xtag}{\@ifstar{\@tag}{\tag@}}
\newcommand{\@tag}{\stepcounter{equation}\tag{\theequation}}
\newcommand{\tag@}[1]{\@tag}
\makeatother
\numberwithin{equation}{section}

\makeatletter
\newcommand{\extraspace}{\@ifstar{\hsp{0cm}\\[-.2cm]}{\hsp{0cm}\\[-.4cm]}}
\newcommand{\xeqref}[2]{
\expandafter\ifx\csname r@#2\endcsname\relax
  \protect\G@refundefinedtrue
  \nfss@text{\reset@font\bfseries (??) R#1}%
  \@latex@warning{%
      Reference `#2' on page \thepage \space undefined%
  }%
\else
  \expandafter\expandafter\expandafter\Hy@setref@link\csname r@#2\endcsname\@empty\@empty\@nil{{\textup{\tagform@{\ref*{#2}}} R#1}\expandafter\@gobble\@firstoffive}
  \fi}
\newcommand{\EXP}{\qg}
\newcommand{\xfrl}{\@ifstar{\frl*\EXP}{\frl\EXP}}
\makeatother

\lonecnt{probl}

\def\L{\mathbb L}

\title{Weak and dissipative solutions for the Hasegawa-Mima equation}
\author{Michele Gorini\fn{SNS - Scuola Normale Superiore, Via dei Consoli del Mare, 2, 56126, Pisa, Italy; \\ \textcolor{white}{abc} E-mail: michele.gorini@sns.it}}
\date{\vsp{-\baselineskip}}

\begin{document}
\maketitle

\begin{abstract}
We consider the Hasegawa-Mima equation in its ``Euler-like'' velocity form:
\[\pd_t(u-\Dg^{-1}u)+(u\per\grad)u-u^\perp\log n_0=0,\]
$n_0$ being the time-independent function appearing in the particle count $n=n_0e^{\xfr{e\fg}{T}}$, and $u$ being the drift velocity $\grad^\perp\fg=-\grad\fg\x\hat z$. Adapting the notion from Lions' book on the Euler equations, we prove the existence of dissipative solutions for this equation for any $L^2$ divergence free initial condition $w\in L^2(D)$, for $D=\T^2$ and $D\sbs\R^2$ a bounded $\s{C}^1$ domain.
\end{abstract}

\sect{Introduction}
The Hasegawa-Mima (HM) equation is a fundamental model for drift-wave turbulence in strongly magnetized plasmas. Introduced in the seminal work of Hasegawa and Mima \cite{HM}, it describes the nonlinear evolution of electrostatic potential fluctuations in the plane perpendicular to a strong magnetic field under the assumption of adiabatic electron response. In nondimensional form in two dimensions, the potential form of the equation reads:
\begin{equation}
\pd_t(\Dg\fg-\fg)+(\grad^\perp\fg\per\grad)(\Dg\fg-\log n_0)=0, \label{eq:HM}
\end{equation}
where $\fg$ denotes the electrostatic potential, the electron count is $n=n_0(x,y)e^{\xfr{e\fg}{T_e}}$, and $\grad^\perp\fg\coloneq(-\pd_2\fg,\pd_1\fg)$.
\nipar Since we are envisioning the poloidal section of a tokamak with a toroidal magnetic field $B$ of unitary magnitude $|B|=1$, $\grad^\perp\fg$ is the incompressible drift velocity $E\x B=-\grad\fg\x\hat z$. The equation can be recast in terms of $u\coloneq\grad^\perp\fg$, i.e. in its velocity form:
\begin{equation}
\pd_t(u-\Dg^{-1}u)+(u\per\grad)u-u^\perp\log n_0=0, \label{eq:ELHM}
\end{equation}
where $u^\perp=(-u_2,u_1)$.
\nipar Since $\pd_tn_0=0=(\grad^\perp\fg\per\grad)\fg$, equation \eqref{eq:HM} can be rewritten as an active scalar transport equation:
\[\pd_tq+(u\per\grad)q=0,\]
where $u$ is the velocity introduced above and $q\coloneq\Dg\fg-\fg-\log n_0$. It is thus closely related to the two-dimensional Euler equations in vorticity form and to quasi-geostrophic shallow water models. If $n_0$ is assumed to be constant, we can reduce the vorticity to $q'\coloneq\Dg\fg-\fg$, and thus, since the elliptic operator $-(I-\Dg)^{-1}$ is smoothing of order two, the velocity field $u$ is one derivative more regular than the transported scalar $q'$. Even without this reduction, formal conservation of energy holds, similarly to the Euler case:
\[\Der{t}\xints D{}\abs{\fg}^2+\abs{\grad\fg}^2\diff x=0\qquad\qquad\Der{t}\xints D{}\abs{(-\Dg)^{-\xfr12}u}^2+\abs{u}^2\diff x=0.\]
Global existence of weak solutions in two space dimensions which propagate the $L^2$ regularity of $(I-\Dg)\fg$ has been established in periodic settings via Galërkin approximations and compactness arguments, e.g. in \cite{KN}, which however does not address conservation of energy and assumes $n_0=n_0(x)$ is independent of $y$. In bounded domains $D\sbs\R^2$, we find \cite{FT}, where finite-time solutions are proved to exist assuming $\log n_0\in L^p(D)$ for some $p$, with uniqueness only in some cases. Those solutions satisfy a bound of the form $\|\fg-\Dg\fg\|_p\lsim\|\fg-\Dg\fg\|_0+\|\log n_0\|_p$. In $\R^2$, finite-time existence and uniqueness was obtained in \cite{GH} with the same assumption on $\log n_0$ as in \cite{KN}.
\nipar As in the case of the incompressible Euler equations, the nonlinear structure of the equation raises compactness issues. In particular, weak limits of approximate solutions (for instance, vanishing viscosity regularizations or spectral truncations) need not satisfy the equation in the classical weak sense due to possible oscillation defects in the quadratic term $(\grad^\perp\fg\per\grad)\Dg\fg$.
\nipar For the Euler equations, this difficulty led Lions to introduce the notion of dissipative solutions in \cite{Lions}, defined through a relative energy inequality with respect to smooth test flows. This concept provides a framework that is stable under weak convergence, allows for possible energy defects, and satisfies weak-strong uniqueness: as long as a classical solution exists, it is unique within the dissipative class. Dissipative and, more generally, measure-valued solutions have since played a central role in the modern analysis of fluid equations, particularly in connection with vanishing viscosity limits and turbulence.
\nipar The aim of the present paper is to develop an analogous theory for the Hasegawa-Mima equation. We introduce notions of dissipative solutions for \eqref{eq:HM} and \eqref{eq:ELHM}, and remark on how one obtains dissipative solutions of one from those of the other. Weak-strong uniqueness for both notions follows directly from the bound that defines the dissipative solutions.
\nipar The main result of this work establishes the global-in-time existence of dissipative solutions of \eqref{eq:ELHM} for arbitrary $L^2$ divergence-free initial data. The construction proceeds via viscous regularization of \eqref{eq:ELHM}, uniform energy estimates, and a compactness argument. Passing to the limit, we obtain a solution satisfying a relative energy inequality with respect to arbitrary smooth divergence-free test functions. As an intermediate step, we prove that the viscous form of \eqref{eq:ELHM} satisfies a result akin to Leray's result \cite{Ler} for the Navier-Stokes equations.
\nipar Beyond its intrinsic analytical interest, the dissipative formulation provides a natural setting for studying singular limits and approximation schemes for the Hasegawa-Mima equation. It is especially relevant in regimes where small-scale oscillations may develop and where classical weak convergence methods fail to identify a limiting equation without defect measures.
\nipar The paper is organized as follows. In Section \ref{ELHMderiv}, we derive the velocity form \eqref{eq:ELHM}. Section \ref{WeakSol} introduces the notions of weak solutions for the two forms, and shows the problems they have. Section \ref{DS} introduces dissipative solutions for both forms (for \eqref{eq:HM} in \ref{DSP} and for \eqref{eq:ELHM} in \ref{DSEL}) and comments on how they relate to each other and the weak-strong uniqueness property they satisfy (\ref{RelWSU}). In Section \ref{BD} we give some remarks on additional work needed to prove our results in bounded domains $D\sbs\R^2$ as opposed to the 2D torus $\T^2$.  In Section \ref{LLR}, we prove a Leray-like result for the viscous form of \eqref{eq:ELHM}, which we introduce there. Section \ref{ExistDiss} uses this result as a starting point to prove the existence of dissipative solutions for \eqref{eq:ELHM}.

\sect{The derivation of velocity-form Hasegawa-Mima}\label{ELHMderiv}
Before proceeding, let us recall a few identities from Appendix \ref{AppA}, where $f_{3D}(x,y,z)\coloneq f(x,y)$ for any two-variable function.
\begin{align}
\grad^\perp\fg={}&-\grad\fg_{3D}\x k \\
(\opn{curl}u)\hat z\coloneq(\pd_1u_2-\pd_2u_1)\hat z={}&\opn{curl}u_{3D} \\
\opn{curl}\opn{curl}u={}&2\grad\opn{div}u-2\Dg u \\
\opn{div}\grad^\perp\fg={}&0 \label{eq:divgp} \\
\opn{curl}\grad\fg={}&0 \label{eq:curlgrad} \\
\opn{div}u={}&\opn{curl}u^\perp \\
\grad^\perp\opn{curl}u={}&\Dg u-\grad\opn{div}u. \label{eq:vattela}
\end{align}
We can also use $\grad\per$ for the divergence and $\grad^\perp\per$ for the curl, which essentially reduces \eqref{eq:divgp} and \eqref{eq:curlgrad} to:
\[\grad\per\grad^\perp=-\grad^\perp\per\grad=0.\]
With this definition of the curl of a 2D field, the immediate recasting of \eqref{eq:HM} in terms of $u\coloneq\grad^\perp\fg$ is as follows:
\[\pd_t(\opn{curl}u-\opn{curl}\Dg^{-1}u)+(u\per\grad)(\opn{curl}u-\log n_0)=0.\]
Aiming to collect a curl operator, we focus on the nonlinearity:
\[\opn{curl}[(u\per\grad)u]=\sum_i\opn{curl}[u_i\pd_iu]=(u\per\grad)\opn{curl}\per u+\sum_i(\grad^\perp u_i)\pd_iu.\]
The second term in this expression vanishes. Indeed:
\begin{align*}
\sum_i\grad^\perp u_i\per\pd_iu={}&-\pd_2u_1\pd_1u_1+\pd_1u_1\pd_1u_2-\pd_2u_2\pd_2u_1+\pd_1u_2\pd_2u_2={} \\
{}={}&\pd_2u_2(\pd_1u_2-\pd_2u_1)+\pd_1u_1(\pd_1u_2-\pd_2u_1)=\opn{curl}u\per\underbrace{\opn{div}u}_0.
\end{align*}
Since $\opn{div}u=\opn{div}\grad^\perp\fg=0$ as seen above, we have:
\[(u\per\grad)\log n_0=\opn{div}(u\log n_0)=\opn{curl}(u^\perp\log n_0).\]
Therefore, our equation reads:
\[\opn{curl}[\pd_t(u-\Dg^{-1}u)+(u\per\grad)u-u^\perp\log n_0]=0.\]
Assuming the argument of the curl vanishes yields \eqref{eq:ELHM}. Coupling this with incompressibility, i.e. $\opn{div}u=0$, implies, by \eqref{eq:vattela}, that $\grad^\perp\opn{curl}\Dg^{-1} u=u$, i.e. $u=\grad^\perp\fg$ for $\fg\coloneq\opn{curl}\Dg^{-1}u$. This means that, if $u$ is incompressible and solves \eqref{eq:ELHM}, then $\fg\coloneq\opn{curl}\Dg^{-1}u$ solves \eqref{eq:HM}.

\sect{Weak solutions to potential-form and velocity-form Hasegawa-Mima}\label{WeakSol}
Assume $\fg\in\Cinf(D)$ solves the potential form \eqref{eq:HM}, and $\qg(x,y)\in\Cinf_c(D)$. We then have:
\begin{equation}
\Der{t}\ang{\fg,\Dg\qg-\qg}-\ang{\grad^\perp\fg\per\grad\qg,\Dg\fg-\log n_0}=0, \label{eq:HMEw}
\end{equation}
where we used $\qg\in\Cinf_c$ to justify the absence of boundary terms when integrating by parts, and the notations:
\[\grad\fg\coloneq(\pd_1\fg,\pd_2\fg)=\sum_1^2\pd_i\fg e_i\qquad\grad^\perp\fg\coloneq(-\pd_2\fg,\pd_1\fg)=\sum_1^2(-1)^i\pd_{(i+1)}\fg e_i,\qquad\pd_3\coloneq\pd_1\]
Indeed, multiplying \eqref{eq:HM} by $\qg$ and integrating in space over $D$, we obtain:
\[\underbrace{\ang{\pd_t(\Dg\fg-\fg),\qg}}_I+\underbrace{\ang{\grad^\perp\fg\per\grad(\Dg\fg-\log n_0),\qg}}_{II}=0.\]
Manipulating term $I$:
\[\ang{\pd_t(\Dg\fg-\fg),\qg}=\Der{t}\ang{\Dg\fg-\fg,\qg}-\ang{\Dg\fg-\fg,\pd_t\qg}=\Der{t}\ang{\fg,\Dg\qg-\qg}-\ang[b]{\fg,\overbrace{\pd_t(\Dg\qg-\qg)}^0}.\]
Manipulating term $II$:
\begin{align*}\ang{\grad^\perp\fg\per\grad(\Dg\fg-\log n_0),\qg}={}&\ang{-\pd_2\fg\per\pd_1(\Dg\fg-\log n_0)+\pd_1\fg\per\pd_2(\Dg\fg-\log n_0),\qg}={} \\
{}={}&\cancel{\ang{\pd_1\pd_2\fg(\Dg\fg-\log n_0),\qg}}+\ang{\pd_2\fg(\Dg\fg-\log n_0),\pd_1\qg}+{} \\
&{}-\cancel{\ang{\pd_2\pd_1\fg(\Dg\fg-\log n_0),\qg}}-\ang{\pd_1\fg(\Dg\fg-\log n_0),\pd_2\qg}={} \\
{}={}&-\ang{\grad^\perp\fg\per\grad\qg,\Dg\fg-\log n_0}.
\end{align*}
This prompts the following definition of weak solutions.
\xbegin{defi}[Weak solution]
Let $D_T$ be a closed time interval or $\R^+$, and $D$ be either $\T^2$, $\R^2$, or an open domain in $\R^2$. A function $\fg\in L^2(D_T;H^2(D)\!)$ is a \emph{weak solution to \eqref{eq:HM}} if it satisfies \eqref{eq:HMEw} for all $\qg\in\Cinf_c(D)$.
\xend{defi}
Choosing a time-dependent $\qg(t,x,y)$, we would obtain:
\begin{equation}
\Der{t}\ang{\fg,\Dg\qg-\qg}-\ang{\fg,\pd_t(\Dg\qg-\qg)}-\ang{\grad^\perp\fg\per\grad\qg,\Dg\fg-\log n_0}=0. \label{eq:HMEt}
\end{equation}
This would prompt a somewhat more restrictive definition, where TDT stands for Time-Dependent Test functions and is introduced for disambiguation purposes.
\xbegin{defi}[TDT-weak solution]
In the same setting as the previous definition, $\fg$ is a \emph{TDT-weak solution to \eqref{eq:HM}} if it satisfies \eqref{eq:HMEt} for all $\qg\in\Cinf(D_T\x D)$.
\xend{defi}
Both of these requires $\fg\in L^2(D_T;H^2(D))$. Indeed, even integrating by parts again in $\ang{\grad^\perp\fg\per\grad\qg,\Dg\fg}$ will not help since:
\[\ang{\grad^\perp\fg\per\grad\qg,\Dg\fg}=-\ang{\grad\fg\per\grad^\perp\qg,\Dg\fg}=\ang{\grad(\grad\fg\per\grad^\perp\qg),\grad\fg},\]
i.e. the integration by parts does not cancel the $\grad^\perp\fg$, leading once again to second-order derivatives of $\fg$. The question regarding the possibility of a definition of $H^1$ weak solutions then comes natural.
\nipar The equivalent of this question for \eqref{eq:ELHM} would be $L^2$ solutions. An approach analogous to the one attempted above for \eqref{eq:HM}, when applied to \eqref{eq:ELHM}, yields:
\begin{align}
0={}&\ang{\pd_t(u-\Dg^{-1}u),v}+\ang{(u\per\grad)u,v}-\ang{u^\perp\log n_0,v}={} \notag \\
{}={}&\Der{t}\ang{u-\Dg^{-1}u,v}-\ang{u-\Dg^{-1}u,\pd_tv}+{} \notag \\
&{}-\ang[B]{\overbrace{\opn{div}u}^0,u\per v}-\ang{u,(u\per\grad)v}+\ang{u,v^\perp\log n_0}. \label{eq:ELHMw}
\end{align}
Integrating in time:
\begin{multline}
\ang{u-\Dg^{-1}u,v}(T)-\ang{u-\Dg^{-1}u,v}(0)+{} \\
-\xints0T\ang{u-\Dg^{-1}u,\pd_tv}\diff t-\xints0T\ang{u,(u\per\grad)v}\diff t+\xints0T\ang{u,v^\perp\log n_0}\diff t. \label{eq:ELHMwint}
\end{multline}
This would then prompt the following definition.
\xbegin{defi}[Weak solution of velocity-form Hasegawa-Mima]
With $D$ defined as in the previous definitions, $u\in\s{C}^0(D_T;L^2_w(D)\!)$ is a \emph{weak solution to \eqref{eq:ELHM}} if it satisfies \eqref{eq:ELHMwint} for all $v\in\Cinf_c(D)$.
\xend{defi}
This only requires $L^2(D)$, but it is not stable under weak convergence. Indeed, if $u_n\rightharpoonup u$ in $L^2$ uniformly in time, the nonlinear term may not pass to the limit, whereas all other terms will, assuming that $\log n_0$ is not too badly behaved ($\log n_0\in L^2$ is enough). To guarantee this passage, we would need to assume uniform-in-time weak convergence in $L^4$ and $L^2$ (i.e. convergence in $\s{C}^0_t(L^2_w)_x\cap\s{C}^0_t(L^4_w)_x$). It would therefore be of interest to develop a different notion of solution which requires less regularity for this kind of stability.

\sect{Dissipative solutions: definitions, equivalence, and weak-strong uniqueness}\label{DS}
\ssect{Dissipative solutions to potential-form Hasegawa-Mima}\label{DSP}
Suppose $\fg,\yg,f\in\Cinf$ satisfy:
\begin{align*}
\pd_t(\Dg\fg-\fg)+\grad^\perp\fg\per\grad(\Dg\fg-\log n_0)={}&0 \\
\pd_t(\Dg\yg-\yg)+\grad^\perp\yg\per\grad(\Dg\yg-\log n_0)={}&f_{\yg}.
\end{align*}
For our convenience, we rearrange these equations as follows:
\begin{align}
\pd_t(\Dg\fg-\fg)={}&-\grad^\perp\fg\per\grad(\Dg\fg-\log n_0) \label{eq:eqphi} \\
\pd_t(\yg-\Dg\yg)={}&\grad^\perp\yg\per\grad(\Dg\yg-\log n_0)-f_{\yg}, \label{eq:eqpsi}
\end{align}
Noting that:
\begin{align}
\xfr12\Der{t}\abs{f}^2={}&\xfr12\Der{t}f^2=f\pd_tf \label{eq:Dersq} \\
\xfr12\xints D{}\Der{t}\abs{\grad f}^2\diff x={}&\xints D{}\grad f\pd_t\grad f\diff x=-\xints D{}f\pd_t\Dg f\diff x+\xints{\mcl{\pd D}}{}f\pd_t\grad f\per\hg\diff\sg, \label{eq:Dersqgr}
\end{align}
we conclude that, if $(\fg,\yg)|_{\pd D}\equiv0$:
\begin{align*}
\xfr12\Der{t}\xints D{}\2\pa{\abs{\grad(\fg-\yg)}^2+\abs{\fg-\yg}^2}\diff x\pux*{\eqref{eq:Dersq}+\eqref{eq:Dersqgr}}{=}{}&-\xints D{}(\fg-\yg)\pd_t(\Dg(\fg-\yg)-(\fg-\yg)\!)\diff x={} \\
{}={}&-\xints D{}(\fg-\yg)[\pd_t(\Dg\fg-\fg)+\pd_t(\yg-\Dg\yg)\!)]\diff x={} \\
{}\pux*[\big]{\eqref{eq:eqphi}+\eqref{eq:eqpsi}}{=}{}&-\xints D{}(\fg-\yg)\sq{-\grad^\perp\fg\per\grad(\lbar{\Dg\fg}-\ubar{\log n_0})}\diff x+{} \\
&{}-\xints D{}(\fg-\yg)\sq{\grad^\perp\yg\per\grad(\lbar{\Dg\yg}-\ubar{\log n_0})-\ubar{f_{\yg}}}\diff x={} \\
{}={}&\underbrace{{}-\xints D{}(\fg-\yg)\sq{\ubar{\grad^\perp(\fg-\yg)\per\grad\log n_0-f_{\yg}}}\diff x}_I+{} \\
&{}+\underbrace{\xints D{}(\fg-\yg)\sq{\grad^\perp\yg\per\grad\Dg(\fg-\lbar\yg)}\diff x}_{II}+{} \\
&{}+\underbrace{\xints D{}(\fg-\yg)\sq{\grad^\perp(\lbar{\fg}-\yg)\per\grad\Dg\fg}\diff x}_{III}.
\end{align*}
Recall the following two identities from Appendix \ref{AppA}:
\begin{align}
\grad f\per\grad^\perp f=-\pd_2f\per\pd_1f+\pd_1f\per\pd_2f=0 \label{eq:GGP} \\
\opn{div}\grad^\perp f=\opn{div}\pa{\e[3]\mat{c}-\pd_2f \\ \pd_1f\emat\e[3]}=-\pd_1\pd_2f+\pd_2\pd_1f=0. \label{eq:DGP}
\end{align}
Thus, integrating by parts, we see that:
\begin{align*}
III={}&\xints D{}(\fg-\yg)\sq{\grad^\perp(\fg-\yg)\per\grad\Dg\fg}\diff x={} \\
{}\pux[\big]{IBP}{=}{}&-\xints D{}\underbrace{\grad(\fg-\yg)\per\grad^\perp(\fg-\yg)}_{\text{0 due to \eqref{eq:GGP}}}\Dg\fg\diff x-\xints D{}(\fg-\yg)\underbrace{\opn{div}\grad^\perp(\fg-\yg)}_{\text{0 due to \eqref{eq:DGP}}}\Dg\fg\diff x=0.
\end{align*}
Integrating $II$ by parts twice:
\begin{align*}
II\pux*[\big]{IBP+\eqref{eq:DGP}}{=}{}&-\xints D{}\!\grad(\fg-\yg)\grad^\perp\yg\Dg(\fg-\yg)\diff x={} \\
{}\pux*[\big]{IBP}{=}{}&\xints D{}\!\grad[\grad(\fg-\yg)\grad^\perp\yg]\grad(\fg-\yg)\diff x-\xints{\pd D}{}\e[3]\sq{\grad(\fg-\yg)\per\grad^\perp\yg}\grad(\fg-\yg)\per\hg\diff\sg={} \\
{}={}&-\xints{\pd D}{}\e[3]\sq{\grad(\fg-\yg)\per\grad^\perp\yg}\grad(\fg-\yg)\per\hg\diff\sg+{} \\
&{}+\sum_1^2\xints D{}\!\grad\pd_i(\fg-\yg)\grad^\perp\yg\pd_i(\fg-\yg)\diff x+\sum_1^2\xints D{}\!\grad(\fg-\yg)\grad^\perp\pd_i\yg\pd_i(\fg-\yg)\diff x\eqcolon{} \\
{}\eqcolon{}&B+A_1+A_2+\sum_1^2\xints D{}\grad(\fg-\yg)\grad^\perp\pd_i\yg\pd_i(\fg-\yg)\diff x.
\end{align*}
Noting that $A_1=0=A_2$:
\begin{align*}
A_i\coloneq{}&\xints D{}\!\grad\pd_i(\fg-\yg)\grad^\perp\yg\pd_i(\fg-\yg)\diff x\pux*[\big]{IBP}{=}-\xints D{}\pd_i(\fg-\yg)\opn{div}[\grad^\perp\yg\pd_i(\fg-\yg)]\diff x={} \\
{}={}&-\xints D{}\pd_i\underbrace{\opn{div}\grad^\perp\yg}_{\text{0 due to \eqref{eq:DGP}}}\pd_i(\fg-\yg)\diff x-\xints D{}\pd_i(\fg-\yg)\grad^\perp\yg\grad\pd_i(\fg-\yg)\diff x=-A_i.
\end{align*}
and putting everything together, we conclude that, if $\fg,\yg\in\Cinf$ satisfy \eqref{eq:eqphi} and \eqref{eq:eqpsi}, and $(\fg-\yg)|_{\pd D}\equiv0$, we have:
\begin{align*}
\xfr12\Der{t}\xints D{}\2\pa{\abs{\grad(\fg-\yg)}^2+\abs{\fg-\yg}^2}\diff x={}&{}-\xints D{}(\fg-\yg)\sq{\grad^\perp(\fg-\yg)\per\grad\log n_0-f_{\yg}}\diff x+{} \\
&{}+\sum_1^2\xints D{}\!\grad(\fg-\yg)\grad^\perp\pd_i\yg\pd_i(\fg-\yg)\diff x+{} \\
&{}-\xints{\pd D}{}\e[3]\sq{\grad(\fg-\yg)\per\grad^\perp\yg}\grad(\fg-\yg)\per\hg\diff\sg.
\end{align*}
We thus conclude the estimate:
\begin{align*}
\Der{t}\xints D{}\2\pa{\abs{\grad(\fg-\yg)}^2+\abs{\fg-\yg}^2}\diff x\leq{}&2\per\!\sum_1^2\xints D{}\!\abs{\grad(\fg-\yg)}\per\e[5]\overbrace{\abs{\grad^\perp\pd_i\yg}}^{\leq\norm{\grad^\perp\pd_i\yg}_{L^\8}}\e[5]\per\overbrace{\abs{\pd_i(\fg-\yg)}}^{\leq\abs{\grad(\fg-\yg)}}\diff x+{} \\
&{}+2\per\e\xints D{}\!\abs{\fg-\yg}\sq[x]{\underbrace{\abs{\grad^\perp(\fg-\yg)}}_{=\abs{\grad(\fg-\yg)}}\per\e[4]\underbrace{\abs{\grad\log n_0}}_{\leq\norm{\grad\log n_0}_{L^\8}}}\diff x+2\per\e\xints D{}(\fg-\yg)f_{\yg}\diff x+{} \\
&{}+\xints{\mcl{\pd D}}{}\!\abs{\grad(\fg-\yg)}\underbrace{\abs{\grad^\perp\yg}}_{{}=\abs{\grad\yg}}\abs{\grad(\fg-\yg)}\diff\sg\leq{} \\
\intertext{}
{}\pux{$ab\leq\xfr12(a^2+b^2)$}{\leq}{}&2\per\!\sum_1^2\norm{\grad^\perp\pd_i\yg}_{L^\8}\per\xints D{}\abs{\grad(\fg-\yg)}^2\diff x+2\per\e\xints D{}(\fg-\yg)f_{\yg}\diff x+{} \\
&{}+\norm{\grad\log n_0}_{L^\8}\pa{\xints D{}\abs{\fg-\yg}^2\diff x+\xints D{}\abs{\grad(\fg-\yg)}^2\diff x}+{} \\
&{}+\norm{\grad^\perp\yg}_{L^\8}\abs{\pd D}\pa{\xints D{}\abs{\fg-\yg}^2\diff x+\xints D{}\abs{\grad(\fg-\yg)}^2\diff x}\leq{} \\
{}\leq{}&g_{\yg}(t)\per\xints D{}\abs{\grad(\fg-\yg)}^2+\abs{\fg-\yg}^2\diff x+2\per\e\xints D{}(\fg-\yg)f_{\yg}\diff x,
\end{align*}
where:
\[g_{\yg}(t)\coloneq2\per\!\sum_1^2\|\grad^\perp\pd_i\yg\|_{L^\8}+2\|\grad\log n_0\|_{L^\8}+\norm{\grad^\perp\yg}_{L^\8}\abs{\pd D}.\]
To simplify the equations, we introduce:
\[K_{\fg,\yg}(t)\coloneq\xints D{}\abs{\grad(\fg-\yg)}^2+\abs{\fg-\yg}^2\diff x\qquad F_{\yg}(t)\coloneq2\per\e\xints D{}(\fg-\yg)f_{\yg}\diff x.\]
We can then write the above estimate as:
\begin{equation}
K_{\fg,\yg}'(t)\leq g_{\yg}(t)K_{\fg,\yg}(t)+F_{\yg}(t). \label{eq:D1}
\end{equation}
Multiplying by $e^{I_a^t(g_{\yg})}$, where $I_a^t(g)=\int_a^tg(t)\diff t$, we have:
\[e^{-I_a^t(g_{\yg})}K_{\fg,\yg}'-e^{-I_a^t(g_{\yg})}g_{\yg} K_{\fg,\yg}\leq e^{-I_a^t(g_{\yg})}F_{\yg}.\]
Integrating in time:
\[e^{-I_a^t(g_{\yg})}K(t)-\overbrace{e^{-I_a^t(g_{\yg})}}^1K(a)\leq\xints ate^{-I_a^s(g_{\yg})}F_{\yg}(s)\diff s,\]
or in other words:
\begin{align*}
K_{\fg,\yg}(t)\leq{}&e^{I_a^t(g_{\yg})}K_{\fg,\yg}(a)+e^{I_a^t(g_{\yg})}\per\xints at\xints D{}e^{-I_a^t(g_{\yg})}(\fg-\yg)(s,x)f_{\yg}(s,x)\diff x\diff s={} \\
{}={}&e^{I_a^t(g_{\yg})}\per\!\xints D{}\abs{\grad(\fg-\yg)}^2(a,x)+\abs{\fg-\yg}^2(a,x)\diff x+{} \\
&{}+\xints ate^{I_s^t(g_{\yg})}\per\!\xints D{}(\fg-\yg)(s,x)f_{\yg}(s,x)\diff x\diff s.
\end{align*}
We therefore come to the following definition.
\xbegin{defi}[Dissipative solution of potential-form Hasegawa-Mima]
A \emph{dissipative solution} of \eqref{eq:HM} is a $\fg\in L^\8(D_T;H^1(D)\!)$ such that, for any $\yg\in L^\8(D_T;\s{C}^2(D)\!)$ satisfying \eqref{eq:eqpsi}, setting $g_{\yg}(t)\coloneq2(\|\grad^\perp\pd_1\yg\|_{L^\8}+\|\grad^\perp\pd_2\yg\|_{L^\8}+\|\grad\log n_0\|_{L^\8}+\norm{\grad\yg}_{L^\8}\abs{\pd D})$, the following estimate holds:
\begin{align}
\xints D{}\abs{\grad(\fg-\yg)}^2(t,x)+\abs{\fg-\yg}^2(t,x)\diff x\leq{}&e^{\int_a^tg_{\yg}\diff s}\per\!\xints D{}\abs{\grad(\fg-\yg)}^2(a,x)+\abs{\fg-\yg}^2(a,x)\diff x+{} \notag\\
&{}+\xints ate^{\int_s^tg_{\yg}\diff s}\per\!\xints D{}(\fg-\yg)(s,x)f_{\yg}(s,x)\diff x\diff s. \label{eq:HMdiss}
\end{align}
\xend{defi}
This definition is similar to the one given for the Euler equations in \cite[Section 4.4]{Lions}, and to prove an existence result one would attempt to imitate the proof of the existence result therein (Proposition 4.2). It seems however logical that the more Euler-like velocity-form recasting \eqref{eq:ELHM} will make adapting the strategy used in \cite{Lions} more straightforward.

\ssect{Dissipative solutions to velocity-form Hasegawa-Mima}\label{DissELHM}\label{DSEL}
Proceeding as in the previous subsection, we assume:
\[\begin{sistema}
\pd_t(u-\Dg^{-1}u)+(u\per\grad)u-u^\perp\log n_0=0 \\
\pd_t(v-\Dg^{-1}v)+(v\per\grad)v-v^\perp\log n_0=f_v
\end{sistema},\]
or in a form more convenient for future calculations:
\begin{equation}
\begin{sistema}
\pd_t(u-\Dg^{-1}u)=-(u\per\grad)u+u^\perp\log n_0 \\
-\pd_t(v-\Dg^{-1}v)=(v\per\grad)v-v^\perp\log n_0-f_v
\end{sistema}. \label{eq:2ELHM}
\end{equation}
We aim to find an equivalence between the two notions of dissipative solutions, so we should aim to estimate something equivalent to $\|\fg-\yg\|_{H^1}^2$. Let us then write that in terms of $u,v$, using the integration by parts rules of Appendix \ref{AppA}:
\begin{align*}
\xints D{}\abs{\fg-\yg}^2+\abs{\grad\fg-\grad\yg}^2\diff x={}&\xints D{}\!\abs{\opn{curl}\Dg^{-1}(u-v)}^2+\abs{\grad\opn{curl}\Dg^{-1}(u-v)}^2\diff x={} \\
{}={}&-\xints D{}\!\Dg^{-1}(u-v)\per\sq[B]{-\grad\underbrace{\opn{div}\Dg^{-1}(u-v)}_0+(u-v)}\diff x+{} \\
&{}-\xints D{}\!\opn{curl}\Dg^{-1}(u-v)\per\Dg\opn{curl}\Dg^{-1}(u-v)\diff x+{} \\
&{}+\xints{\mcl{\pd D}}{}\!\opn{curl}\Dg^{-1}(u-v)\Dg^{-1}(u-v)\per\hg^\perp\diff\sg+{} \\
&{}+\xints{\mcl{\pd D}}{}\!\opn{curl}\Dg^{-1}(u-v)\grad\opn{curl}\Dg^{-1}(u-v)\per\hg\diff\sg={} \\
{}={}&\xints D{}(-\Dg)^{-1}(u-v)\per(u-v)-\opn{curl}(u-v)\per\opn{curl}\Dg^{-1}(u-v)\diff x+{} \\
&{}+\xints{\mcl{\pd D}}{}\!\opn{curl}\Dg^{-1}(u-v)[\grad^\perp\opn{curl}\Dg^{-1}(u-v)-\Dg^{-1}(u^\perp-v^\perp)]\per\hg\diff\sg={} \\
{}={}&\xints D{}\!\abs{u-v}^2-(u-v)\grad\underbrace{\opn{div}\Dg^{-1}(u-v)}_0+\abs{(-\Dg)^{-\xfr12}(u-v)}^2\diff x+{} \\
&{}+\xints{\mcl{\pd D}}{}\!\opn{curl}\Dg^{-1}(u-v)[\Dg^{-1}(u-v)-\cancel{(u-v)+(u-v)}]\per\hg^\perp\diff\sg.
\end{align*}
From now on, we assume:
\bi
\item $u|_{\pd D}=0=v|_{\pd D}$ and $\Dg$ is the Dirichlet laplacian -- hence, the boundary term vanishes;
\item $(-\Dg)^{\ag}$ can be defined for non-integer $\ag$ in such a way that $\ang{(-\Dg)^{\ag} f,(-\Dg)^{\bg} g}=\ang{(-\Dg)^{\ag+\bg}f,g}=\ang{f,(-\Dg)^{\ag+\bg}g}$ for any $f,g,\ag,\bg$; on $\T^2$ this can be done by assigning $(-\Dg)^{\ag}$ the Fourier symbol $|\jg|^{2\ag}$, and we will see how this is possible in bounded domains in the next section.
\ei
Under these assumptions, using our rewritten equations \eqref{eq:2ELHM}:
\begin{align*}
\xfr12\Der{t}\xints D{}\abs{u-v}^2+\abs{(-\Dg)^{-\xfr12}(u-v)}^2\diff x={}&\xints D{}(u-v)[\pd_t(u-v)-\pd_t\Dg^{-1}(u-v)]\diff x={} \\
{}={}&\xints D{}(u-v)\pa{(v\per\grad)v-(u\per\grad)u}\diff x+{} \\
&{}+\xints D{}(u-v)\pa{u^\perp\log n_0-v^\perp\log n_0-f_v}\diff x={} \\[-.7em]
\intertext{}
{}={}&\xints D{}\overbrace{(u-v)(u^\perp-v^\perp)}^0\log n_0-(u-v)f_v\diff x+{} \\
&{}-\underbrace{\xints D{}(u-v)(u\per\grad)(u-v)\diff x}_{\eqcolon I}+-\underbrace{\xints D{}(u-v)[(u-v)\per\grad]v\diff x}_{\eqcolon II}.
\end{align*}
Manipulating term $I$:
\[I\coloneq-\xints D{}(u-v)(u\per\grad)(u-v)\diff x=\xints D{}\pd_i[u_i(u-v)_j](u-v)_j\diff x+\e\underbrace{\xints{\mcl{\pd D}}{}\!\abs{u-v}^2u\per\hg\diff\sg}_{0\text{ since }u|_{\pd D}\equiv0}=-I,\]
because $\pd_iu_i=\opn{div}\grad^\perp\fg=0$, and $u_i\pd_i(u-v)_j(u-v)_j=(u-v)(u\per\grad)(u-v)$. Thus:
\[\xfr12\Der{t}\xints D{}\!\abs{u-v}^2+\abs{(-\Dg)^{-\xfr12}(u-v)}^2\diff x=-\xints D{}(u-v)[(u-v)\per\grad]v+(u-v)f_v\diff x.\]
This implies, adding the term with the underbrace to have the same norm of $u-v$ appear on both sides:
\begin{align*}
\xfr12\Der{t}\xints D{}\!\abs{u-v}^2+\abs{(-\Dg)^{-\xfr12}(u-v)}^2\diff x\leq{}&2\norm{\grad v}_{L^\8}\per\!\xints D{}\abs{u-v}^2+\underbrace{\abs{(-\Dg)^{-\xfr12}(u-v)}^2}\diff x+{} \\
&\overbrace{{}-2\xints D{}(u-v)f_v\diff x}^{{}\eqcolon f_{u,v}(t)},
\end{align*}
and thus, calling $\|f\|_{H_{\dg} L}^2\coloneq\|f\|_{L^2}^2+\|(-\Dg)^{-\xfr12}f\|_{L^2}^2$:
\begin{equation}
\norm{u-v}_{H_{\dg} L}^2(t)\leq e^{\xints 0t2\|\grad v\|_\8(s)\diff s}\norm{u_0-v_0}_{H_{\dg} L}^2+\xints0te^{\xints st2\|\grad v\|_\8(\tg)\diff\tg}f_{u,v}(s)\diff s. \label{eq:ELHMdiss}
\end{equation}
We therefore come to the following definition.
\xbegin{defi}[Dissipative solutions to ELHM]
$u\in L^2(D_T;L^2(D)\!)$ is a \emph{dissipative solution to \eqref{eq:ELHM}} if, for all divergence-free $v\in\s{C}^1$, the bound \eqref{eq:ELHMdiss} holds.
\xend{defi}
\clearpage

\ssect{Dissipative solutions of one form of Hasegawa-Mima from those of the other one, and weak-strong uniqueness for both forms}\label{RelWSU}
\xbegin{propo}[Dissipative solutions of both forms]
If $u,v$ satisfy this bound and $u=\grad^\perp\fg,v=\grad^\perp\yg$, then $\fg,\yg$ satisfy the bound for dissipative solutions of \eqref{eq:HM} for $f_{\yg}\coloneq\grad^\perp\per f_v$.
\xend{propo}
\begin{qeddim}
We have already seen that $\|\fg-\yg\|_{H^1}^2=\|u-v\|_{H_{\dg} L}^2$, reducing the above bound to:
\[\norm{\fg-\yg}_{H^1}^2(t)\leq e^{\xints 0t2\|\grad\grad^\perp\yg\|_\8(s)\diff s}\norm{\fg_0-\yg_0}_{H^1}^2+\xints0te^{\xints st2\|\grad\grad^\perp\yg\|_\8(\tg)\diff\tg}f_{u,v}(s)\diff s.\]
The $\|\grad\grad^\perp\yg\|_\8$ in the exponentials was a shorthand $\sup_i\|\pd_i(\grad^\perp\yg)\|_{L^\8}\leq\sum_i\|\pd_i(\grad^\perp\yg)\|_{L^\8}$, the latter being what we have in the bound \eqref{eq:HMdiss}. Finally:
\[f_{u,v}=-2\!\xints D{}(u-v)f_v\diff x=-2\!\xints D{}\!\grad^\perp(\fg-\yg)f_v\diff x=2\!\xints D{}(\fg-\yg)\per\grad^\perp\per f_v\diff x-\xints{\mcl{\pd D}}{}(\fg-\yg)f_v\per\hg^\perp\diff\sg.\]
Since $\fg|_{\pd D}\equiv0\equiv\yg|_{\pd D}$, the boundary term vanishes, which means $f_{u,v}(t)$ is exactly what we had in the bound \eqref{eq:HMdiss}.
\end{qeddim}
Continuing, let us state a trivial weak-strong uniqueness result.
\xbegin{propo}[Weak-strong uniqueness]
If $\fg,u$ are dissipative solutions of \eqref{eq:HM} and \eqref{eq:ELHM} respectively, and $\yg,v$ are respectively a $\s{C}^2$ solution of \eqref{eq:HM} with $\yg(\per,0)=\fg(\per,0)$ and a $\s{C}^1$ solution of \eqref{eq:ELHM} with $v(\per,0)=u(\per,0)$, then $\fg=\yg$ and $u=v$.
\xend{propo}
\begin{qeddim}
Testing $\fg$ against $\yg$ in \eqref{eq:HMdiss}, we see the first term on the RHS is zero because of the identity of the initial data, while the second one is zero since $f_{\yg}=0$, thus implying $\|\fg-\yg\|_{H^1}^2\leq0$ for all times. Similarly, testing $u$ against $v$ in \eqref{eq:ELHMdiss}, since $f_v=0$, we obtain $u=v$.
\end{qeddim}
\vsp{-\baselineskip}
\xbegin{oss}[Why ``dissipative''?]
We use the term \emph{dissipative} mostly because this is the terminology used by Lions, but the term also makes sense on its own, at least for the velocity form. Indeed, testing with $v=0$, the bound \eqref{eq:ELHMdiss} reduces to:
\[\norm{u}_{H_{\dg} L}^2\leq1\per\norm{u_0}_{H_{\dg} L}^2+\xints 0t\!1\per0\,\diff s,\]
which is an energy inequality akin to the one obtained for the Euler equations. If we look at \eqref{eq:HMdiss} for $\yg=0$, the presence of the $\log n_0$ term in the exponentials makes it so the bound reduces to:
\[\norm{\fg}_{H^1}^2\leq F\per\norm{\fg}_{H^1}^2,\]
where $F$ could be greater than 1, thus not quite bounding the energy from above with the initial energy, but almost.
\xend{oss}

\sect{Bounded 2D domains}\label{BD}
The aim of this section is to make sure that the following points do not present problems when working on bounded domains $D$ rather than on the torus $\T^2$:
\ben
\item Sobolev embeddings;
\item An analogue of $e^{ik\per x}$, i.e. an orthonormal basis consisting of eigenfunctions of $\Dg$;
\item Defining fractional laplacians.
\een
The Sobolev embeddings we will need all hold for bounded domains $D$ with $\s{C}^1$ boundary, so we will be assuming such regularity from now on.
\nipar Coming to the existence of an orthonormal basis of $L^2$ consisting of eigenfunctions of $-\Dg$ with zero boundary conditions, i.e. $\{\fg_i\}$ such that $-\Dg\fg_i=\lg_i\fg_i$, $\int_D\fg_i\fg_j\diff x=\dg_{ij}$, and $\fg_i|_{\pd D}\equiv0$, obtaining it requires us to prove that $-\Dg$ is a self-adjoint operator with at least one compact resolvent operator. Self-adjointness is the easy part:
\[\xints D{}\!f\Dg g\diff x=\xints{\mcl{\pd D}}{}\!f\grad g\per\ng-\xints D{}\!\grad f\per\grad g\diff x=\xints{\mcl{\pd D}}{}\!f\grad g\per\ng-\xints{\mcl{\pd D}}{}\!g\grad f\per\ng+\xints D{}\!\Dg f\per g\diff x,\]
and Dirichlet boundary conditions, i.e. $f|_{\pd D}=0=g|_{\pd D}$, ensure that the boundary terms vanish. Hence, $\Dg$ (and thus $-\Dg$) is indeed self-adjoint.
\nipar As for the compactness of its resolvent, if $-\Dg-\lg\opn{Id}$ is injective, then its inverse has image contained in $H^2$. This means that, when we compose with $\opn{id}:H^2\hra L^2$ which is compact, the result is a compact operator by the Rellich-Kondrachov theorem, at least if $D$ has $\s{C}^1$ boundary and this inverse is bounded. In particular, for $\lg=0$, $(-\Dg)^{-1}$ is injective because if $(-\Dg)^{-1}u=(-\Dg)^{-1}v$ then $u=-\Dg(-\Dg)^{-1}u=-\Dg(-\Dg)^{-1}v=v$. Moreover, $\|(-\Dg)^{-1}u\|_{H^2}\lsim\|\Dg(-\Dg)^{-1}u\|_{L^2}=\|u\|_{L^2}$ by the Poincaré inequality applied twice. Thus, $(-\Dg)^{-1}:L^2\to L^2$ is compact.
\nipar Taking our orthonormal basis $\{\fg_i\}_{i\in*\N}$ such that $\int\fg_i\fg_j\diff x=\dg_{ij}$ and any function can be represented as $u=\sum\int u\fg_i\diff x\fg_i$, and moreover $-\Dg\fg_i=\lg_i\fg_i$ for some $\lg_i>0$, then we simply define:
\[(-\Dg)^{\ag} f\coloneq\sum_{i=0}^\8\lg_i^{\ag}\fg_i\xints D{}\!f(x)\fg_i(x)\diff x.\]
The following calculations, where underlined terms equal $\dg_{ij}$, show that all such powers will be self-adjoint:
\begin{align*}
\xints D{}(-\Dg)^{\ag+\bg}f\per g\diff x={}&\xints D{}\e\pa[B]{\sum_{ij}\lg_i^{\ag+\bg}\fg_i(x)\xints D{}\!f(x')\fg_i(x')\diff x'}\pa[B]{\fg_j(x)\xints D{}\!g(x'')\fg_j(x'')\diff x''}\diff x={} \\
{}={}&\sum_{ij}\lg_i^{\ag+\bg}\xints D{}\xints D{}\!f(x')\fg_i(x')g(x'')\fg_j(x'')\ubar{\xints D{}\fg_i(x)\fg_j(x)\diff x}\diff x'\diff x''={} \\[-.5em]
{}={}&\sum_{ij}\lg_i^{\ag}\lg_j^{\bg}\xints D{}\xints D{}\!f(x')\fg_i(x')g(x'')\fg_j(x'')\ubar{\xints D{}\fg_i(x)\fg_j(x)\diff x}\diff x'\diff x''={} \\
\intertext{}
{}={}&\xints D{}\e\pa[B]{\sum_{ij}\lg_i^{\ag}\fg_i(x)\xints D{}\!f(x')\fg_i(x')\diff x'}\pa[B]{\lg_j^{\bg}\fg_j(x)\xints D{}\!g(x'')\fg_j(x'')\diff x''}\diff x={} \\
{}={}&\xints D{}(-\Dg)^{\ag} f\per(-\Dg)^{\bg} g\diff x.
\end{align*}

\sect{Viscous forms and Leray-like result for the velocity form}\label{LLR}
Thanks to the similarity with Euler, it is natural to give \eqref{eq:ELHM} the viscous form below, with the viscous term modeled on the one in the Navier-Stokes equations:
\begin{equation}
\pd_t(u-\Dg^{-1}u)+(u\per\grad)u-u^\perp\log n_0=\ng\Dg u. \label{eq:VELHM}
\end{equation}
The equivalent form for \eqref{eq:HM} is:
\begin{equation}
\pd_t(\Dg\fg-\fg)+(\grad^\perp\fg\per\grad)(\Dg\fg-\log n_0)=\ng\Dg^2\fg. \label{eq:VHM}
\end{equation}
The Leray result \cite{Ler} for the Navier-Stokes equations concerns solutions which satisfy an energy in equality. In order to obtain an analogous inequality for \eqref{eq:VELHM}, let us multiply it by $u$:
\[0=u\per[\pd_t(u-\Dg^{-1}u)+(u\per\grad)u-u^\perp\log n_0-\ng\Dg u]=\Der{t}\abs u^2-u\per\pd_t\Dg^{-1}u+u\per(u\per\grad)u-\ng u\per\Dg u.\]
Integrating in space:
\begin{align*}
\xfr12\Der{t}\xints D{}\abs{u}^2+|(-\Dg)^{-\xfr12}u|^2\diff x={}&\xfr12\xints D{}\sq{\Der{t}|u|^2-u\pd_t\Dg^{-1}u}\diff x={} \\
{}={}&-\xints D{}u\per(u\per\grad)u\diff x+\ng\xints D{}u\per\Dg u\diff x={} \\
{}={}&\xints D{}u_j\pd_i(u_iu_j)\diff x-\underbrace{\xints{\mcl{\pd D}}{}\abs{u}^2u\per\hg\diff\sg}_{0\text{ since }u|_{\pd D}\equiv0}+{} \\
&{}-\ng\xints D{}\abs{\grad u}^2\diff x+\underbrace{\ng\xints{\mcl{\pd D}}{}(\hg\per\grad)u\per u\diff\sg}_{0\text{ since }u|_{\pd D}\equiv0}={} \\
{}={}&\xints D{}\underbrace{\opn{div}u}_0\abs{u}^2\diff x+\underbrace{\xints D{}u(u\per\grad)u\diff x}_0-\ng\xints D{}\abs{\grad u}^2\diff x,
\end{align*}
where the second term is zero because, as can be seen within the calculations above, $\int_Du(u\per\grad)u\diff x=-\int_Du(u\per\grad)u\diff x$. Summing up, we have obtained:
\[\xfr12\Der{t}\xints D{}\abs{u}^2(t,x)+\abs{(-\Dg)^{-\xfr12}u}^2(t,x)\diff x=-\ng\xints D{}\abs{\grad u}^2\diff x.\]
The desired Leray-like result is thus the following.
\xbegin{teor}[Leray-like result for velocity-form Hasegawa-Mima][thm:teor:LLR]
Assume $\log n_0\in L^{2+\eg}$ for some $\eg$. For any $u_0\in L^2$, there exists at least one global-in-time weak solution $u$ of \eqref{eq:VELHM} such that the energy inequality below holds, where $\|\per\|_{H_{\dg} L}$ is defined as in subsection \ref{DissELHM} and $\|u(t,\per)\|_{\dot H^1}\coloneq\|\grad u(t,\per)\|_{L^2}$:
\begin{equation}
\norm{u(t,\per)}_{H_{\dg} L}^2+\ng\xints0t\!\norm{u(t,\per)}_{\dot H^1}^2\diff s\leq\norm{u_0}_{H_{\dg} L}^2. \label{eq:VELHMadm}
\end{equation}
\xend{teor}
\begin{qeddim}[\xbegin{defi}[Leray solutions for ELHM]
Solutions as in this theorem will be called \emph{Leray solutions to the velocity-form Hasegawa-Mima equation}.
\xend{defi}\vsp{-\baselineskip}]
We proceed by adapting the strategy of \cite[Appendix A]{CDLDR}. For $D=\T^2$, we set $\fg_j\coloneq e^{ij\per x}$, and for $D\sbse\R^2$ bounded we consider the basis from the previous section. We then consider $\P_K:L^2(D)\to L^2(D)$ which truncates the abstract Fourier series:
\[P_K(f)(x)\coloneq\sum_{\mcl{|k|\leq K}}\hat f_k\fg_k,\]
and extend it to vector functions by applying it to each component. Observe that this operator commutes with the derivatives and any power of the laplacian.
\nipar For an arbitrary $w_0\in L^2$, we then consider the reguralized Cauchy problem:
\[\begin{sistema}
\pd_t(w_K-\Dg^{-1}w_K)+\opn{div}P_K(w_K\otimes w_K)-w_K^\perp\log n_0-\ng\Dg w_K=0 \\
\opn{div}w_K=0 \\
w_K|_{\pd D}\equiv0 \\
\xints D{}w_k\diff x\equiv0 \\
w(\per,0)=P_K(w_0)
\end{sistema}.\]
This reduces to a system of ODEs for the Fourier coefficients $(\widehat{w_K})_k$ of the solution $w_K$, ensuring local well-posedness. On the other hand, multiplying by $w_K$ and using Plancherel's theorem, we easily see that:
\[\Der{t}\xints D{}|w_K|^2+|(-\Dg)^{-\xfr12}w_K|^2\diff x=-2\ng\per\!\xints D{}\abs{\grad w_K}^2\diff x,\]
proving therefore that any solution stays bounded in $L^2$ in its interval of existence. A standard continuation argument then yields the existence of a global-in-time solution of the ODE system for the $\hat w_K(t)$, i.e. that the above Cauchy problem is globally solvable. Letting $w_K$ be such solution, we observe that:
\begin{align*}
\xfr12\xints D{}\abs{w_K}^2+\abs{(-\Dg)^{-\xfr12}w_K}^2\diff t+\ng\per\!\xints0t\xints D{}\abs{\grad w}^2\diff x={}&\xfr12\xints D{}\abs{P_K(w_0)}^2(x)+\abs{(-\Dg)^{-\xfr12}P_K(w_0)}^2\diff x\leq{} \\
{}\leq{}&\xfr12\xints D{}\abs{w_0}^2(x)+\abs{(-\Dg)^{-\xfr12}w_0}^2\diff x.
\end{align*}
The sequence $\{w_K\}_{K\in*\N}$ is thus bounded in $L^2\cap H^{-1}(D\x[0,T])$ for any $T>0$, meaning we can extract a subsequence, not relabeled, such that $w_K\rightharpoonup v$ in $L^2\cap H^{-1}$. With a standard diagonal argument, we can then assume that the subsequence converges in $(L^2\cap H^{-1})(D\x[0,T])$ for every $T>0$. We now wish to show this subsequence converges locally strongly in $L^2\cap H^{-1}(D\x[0,T])$, which would show $v$ is a solution.
\xbegin{lemma}[Sobolev embedding]
For $D\sbs\R^2$ bounded or $D=\T^2$, $H^1(D)\hra L^p(D)$ compactly for all $p<\8$.
\xend{lemma}
Since we have a uniform estimate for $w_K$ in $L^2_tH^1_x(\R^+\x D)$, by the embedding lemma above we conclude $w_K$ is bounded in all $L^2_tL^p_x$, and by interpolating with the $L^\8_tL^2_x$ bound we conclude there exists $\gg>2$ such that $w_K$ is bounded in $L^{\gg}_tL^{\gg}_x$.
\nipar Fix $T>0$ and define:
\[A_{K,J}\coloneq\norm{(I-\Dg^{-1})(w_K-w_J)}_{L^2_tL^2_x}.\]
Let $\eg>0$ be given. We want to show that there exists $N_{\eg}\in\N$ such that $A_{K,J}<\eg$ for all $K,J\geq N_{\eg}$. Fix $\fg_{\dg}$ a standard space mollifier and note that:
\[\norm{(I-\Dg^{-1})w_K(\per,t)-(I-\Dg^{-1})w_K\ast\fg_{\dg}(\per,t)}_{L^2}\lsim\dg\norm{(I-\Dg^{-1})w_K(\per,t)}_{H^1}.\]
Now, $w_K$ is bounded in $H^1$. Since $w_K$ is mean-free:
\begin{align*}
\|\Dg^{-1}w_K\|_{H^1}^2={}&\|\Dg^{-1}w_K\|_{L^2}^2+\|\Dg^{-1}Dw_K\|_{L^2}^2\leq\norm{\Dg^{-1}w_K}_{H^2}^2+\norm{\Dg^{-1}Dw_K}_{H^2}^2\lsim{}\\
{}\pux{Poincaré}{\lsim}{}&\norm{D^2\Dg^{-1}w_K}_{L^2}^2+\norm{D^2\Dg^{-1}Dw_K}_{L^2}\lsim\norm{\Dg\Dg^{-1}w_K}_{L^2}^2+\norm{\Dg\Dg^{-1}Dw_K}_{L^2}=\norm{w_K}_{H^1}^2.
\end{align*}
This implies that, for $\dg$ sufficiently small, we have:
\[\norm{(I-\Dg^{-1})w_K\ast\fg_{\dg}-(I-\Dg^{-1})w_K}_{L^2}<\xfr\eg3\qquad\VA K\in\N.\]
Observe now that, by the boundedness of $w_K$ in $L^{\gg}_tL^{\gg}_x$, we have:
\[\norm{P_K(w_K\otimes w_K)}_{L^{\xfr\gg2}_{t,x}}\leq\norm{w_K\otimes w_K}_{L^{\xfr\gg2}_{t,x}}\leq\norm{w_K}_{L^{\gg}_{t,x}}^2\lsim1.\]
Thus, mollifying the equation for $w_K$ we find:
\[\pd_t(I-\Dg^{-1})w_K\ast\fg_{\dg}=-\sum_{i=1}^3f_{i,K}\ast\pd_{x_i}\fg_{\dg}-(w_K^\perp\log n_0)\ast\fg_{\dg}-w_K\ast(-\Dg)\fg_{\dg},\]
where the $f_{i,K}$ are uniformly bounded in $L^{\xfr\gg2}$. Using the estimate $\|\zg\ast\fg_{\dg}\|_{W^{1,\8}}\leq C(\dg)\|\zg\|_{L^{\xfr\gg2}}$ for each time slice, using the embeddings $L^{\gg},L^{2+\eg}\hra L^{\min\{\gg,2+\eg\}}$ for the $\log n_0$ term, we conclude:
\[\xints0T\!\norm{\pd_t(I-\Dg^{-1})w_K\ast\fg_{\dg}(\per,t)}_{W^{1,\8}}^{\xfr{\min\{2+\eg,\gg\}}2}\diff t\leq C(\dg),\]
where $C(\dg)$ is a constant that depends on $\dg$ but not on $K$.
\nipar We can thus regard $t\mapsto(I-\Dg^{-1})w_K\ast\fg_{\dg}(\per,t)$ as a sequence of equicontinuous and equibounded $W^{1,\8}(D)$-valued curves. Let $B_R$ be a (closed) ball in that space so that the images of $w_K\ast\fg_{\dg}$ are all contained inside it. Endowing $B_R$ with the $\|\per\|_\8$ norm, we have a compact metric space, allowing us to apply the Ascoli-Arzelà theorem, concluding the sequence is precompact. Since the limit is $(I-\Dg^{-1})v\ast\fg_{\dg}$(\fn{$w_K\xrightharpoonup{L^2}v$, so $w_K\ast\fg_{\dg}$ can only converge uniformly to $v\ast\fg_{\dg}$, and we have just seen $\Dg^{-1}$ is bounded $H^1\to H^1$ and, by a similar argument, $L^2\to L^2$, so the same argument applies to $\Dg^{-1}w_K\xrightharpoonup{L^2}\Dg^{-1}v$}), we conclude $(I-\Dg^{-1})w_K\ast\fg_{\dg}$ converges uniformly on $D\x[0,T]$. Thus there exists $N$ large enough such that:
\[\norm{(I-\Dg^{-1})w_K\ast\fg_{\dg}-(I-\Dg^{-1})w_J\ast\fg_{\dg}}_{L^2_tL^2_x}<\xfr\eg3.\]
Therefore, for $J,K\geq N$, we have:
\begin{align*}
\norm{(I-\Dg^{-1})w_K-(I-\Dg^{-1})w_J}\leq{}&\norm{(I-\Dg^{-1})w_K-(I-\Dg^{-1})w_K\ast\fg_{\dg}}+{} \\
&{}+\norm{(I-\Dg^{-1})w_K\ast\fg_{\dg}-(I-\Dg^{-1})w_J\ast\fg_{\dg}}+{} \\
&{}+\norm{(I-\Dg^{-1})w_J\ast\fg_{\dg}-(I-\Dg^{-1})w_J}<\eg.
\end{align*}
Hence, $(I-\Dg^{-1})w_K\to(I-\Dg^{-1})v$ strongly in $L^2$. Looking at our abstract Fourier series and using Plancherel's theorem:
\[\norm{w_K}_{L^2}^2=\sum_k\abs{\hat w_{K,k}}^2\leq\sum_k(1+|k|^{-2})^2\abs{\hat w_{K,k}}^2=\norm{(I-\Dg^{-1})w_K}_{L^2}^2,\]
i.e. the convergence of $(I-\Dg^{-1})w_K$ implies that of $w_K$. Since both $w_K$ and $(I-\Dg^{-1})w_K$ converge in $L^2$, we can conclude:
\[\xints D{}\abs{w_K}^2+\abs{(-\Dg)^{-\xfr12}w_K}^2\diff x=\xints D{}w_K(I-\Dg^{-1})w_K\diff x\to\xints D{}v(I-\Dg^{-1})v\diff x=\xints D{}\abs{v}^2+\abs{(-\Dg)^{-1}v}^2\diff x.\]
We can thus conclude that:
\[\xfr12\xints D{}\abs{v}^2+\abs{(-\Dg)^{-\xfr12}v}^2\diff x+\liminf_{K\to\8}\xints 0t\xints D{}\abs{\grad w_K}^2\diff x\diff s\leq\xfr12\xints D{}\abs{\lbar v}^2+\abs{(-\Dg)^{-\xfr12}\lbar v}^2\diff x.\]
As for that $\liminf$, since $w_K$ converges weakly in $L^2H^1$, $Dw_K$ converges weakly in $L^2_tL^2_x$, which yields:
\[\norm{\grad v}_{L^2_tL^2_x}\leq\liminf_K\norm{\grad w_K}_{L^2_tL^2_x}.\]
Thus, the estimate above implies:
\[\xfr12\xints D{}\abs{v}^2+\abs{(-\Dg)^{-\xfr12}v}^2\diff x+\xints 0t\xints D{}\abs{\grad v}^2\diff x\diff s\leq\xfr12\xints D{}\abs{\lbar v}^2\diff x+\abs{(-\Dg)^{-\xfr12}\lbar v}^2,\]
concluding the proof. \vsp{-\baselineskip}
\end{qeddim}
\xbegin{lemma}[Weak continuity][thm:lemma:WCont]
Let $v\in L^\8_t(D_T;L^2_x\cap H^{-1}_x(D)\!)$ and $u\in L^1_{loc}(D\x D_T,\R^{2\x 2})$ be distributional solutions of:
\begin{equation}
\pd_t(v-\Dg^{-1}v)+\opn{div}_xu-v^\perp\log n_0=0. \label{HMbis}
\end{equation}
Then, after redefining $v$ on a set $\tg\sbs D_T$ of measure zero, $v\in\s{C}^0_t(L^2_w)_x$.
\xend{lemma}
\begin{qeddim}[We adapt the strategy of \cite[Appendix A, Lemma 7.1]{DLSz10}. We will then apply this for $u=v\otimes v+\grad v$.]
Consider a countable set $\{\fg_i\}\sbse\Cinf_c(D,\R^2)$ dense in the strong topology of $L^2$. Fix $\fg_i$ and any test function $\xg\in\Cinf_c(D_T)$. Testing \eqref{HMbis} with $\xg(t)\fg_i(x)$ we obtain the following identity:
\begin{equation}
\xints{\mcl{D_T}}{}\Fg_i\pd_t\xg=-\xints{\mcl{D_T}}{}\xg\xints D{}\sq{\ang{u,\grad\fg_i}-v^\perp\fg_i\log n_0}, \label{eq:HMter}
\end{equation}
where $\Fg_i(t)\coloneq\int_D\fg_i(x)\per(v(x,t)-\Dg^{-1}v(x,t))\diff x$. We conclude, therefore, that $\Fg_i'\in L^1$ in the sense of distributions. Hence we can redefine each $\Fg_i$ on a set of times $\tg_i\sbs\pint D_T$ of measure zero in such a way that $\Fg_i$ is continuous. We keep the same notation for these functions, and let $\tg\coloneq\cup_i\tg_i$. Then $\tg\sbs\pint D_T$ is of measure zero and for every $t\in\pint D_T\ssm\tg$ we have:
\[\Fg_i(t)=\xints D{}\fg_i(x)(v(x,t)-\Dg^{-1}v(x,t)\!)\diff x=\xints D{}(I-\Dg^{-1})\fg_i(x)\per v(x,t)\diff x\qquad\VA i.\]
Moreover, with $c\coloneq\|v\|_{L^\8_tL^2_x}$, we have that $|\Fg_i(t)|\leq c\|(I-\Dg^{-1})\fg_i\|_{L^2}\leq Cc\|\fg_i\|_{L^2}$ for all $t\in\pint D_T$, since $\Dg^{-1}:L^2\to L^2$ is compact, as seen above. Therefore, for each $t\in\pint D_T$ there exists a unique bounded linear functional $L_t$ on $L^2(D,\R^2)$ such that $L_t(\fg_i)=\Fg_i(t)$. By the Riesz representation theorem there exists $\bar v(\per,t)\in L^2(D)$ such that:
\bi
\item $\bar v(\per,t)=v(\per,t)$ for all $t\in\pint D_T\ssm\tg$;
\item $\|\bar v(\per,t)\|_{L^2}\leq Cc$ for all $t$;
\item $\int\bar v(x,t)\per(I-\Dg^{-1})\fg_i(x)\diff x=\Fg_i(t)$ for all $t$.
\ei
To conclude, we show that $\bar v\in\s{C}^0(\pint D_T;L^2_w)$, i.e. that $\Fg(t)\coloneq\int_D\bar v(x,t)\fg(x)\diff x$ is continuous on $\pint D_T$ for all $\fg\in L^2(D,\R^2)$. Assume at first that $\fg\in\Cinf$. Then, $(I-\Dg^{-1})^{-1}\fg$ is still smooth, and since the set $\fg_i$ is dense in $L^2(D,\R^2)$, we can find a sequence $\{j_k\}$ such that $\fg_{j_k}\to(I-\Dg^{-1})^{-1}\fg$ strongly in $L^2$, meaning $(I-\Dg^{-1})\fg_{j_k}\to\fg$ strongly by boundedness of $I-\Dg^{-1}:L^2\to L^2$. Then:
\[\abs{\Fg(t)-\Fg_{j_k}(t)}\leq c\norm{\fg_{jk}-\fg}_{L^2}.\]
Therefore $\Fg_{j_k}$ converges uniformly to $\Fg$, from which we derive the continuity of $\Fg$. Since $\Cinf$ is dense in the strong topology, since $\Fg$ is continuous for any $\fg\in\Cinf$, if $\fg\in L^2\ssm\Cinf$, just choose $\{\fg_i\}\sbs\Cinf$ such that $\fg_i\to\fg$ in $L^2$ and argue analogously, concluding $\Fg$ is indeed continuous for all $\fg\in L^2(D,\R^2)$. This shows $V\in\s{C}^0_t(L^2_w)_x$. \placeqed*
\end{qeddim}
\xbegin{cor}[][thm:cor:WCont]
A Leray solution of \eqref{eq:VELHM} is $\s{C}^0_t(L^2_w)_x$ up to redefinition on $\tg\sbs D_T$ of measure 0.
\xend{cor}
Apply the lemma above for $u=v\otimes v+\grad v$, $v$ being the solution of \eqref{eq:VELHM} we are starting from, and the conclusion is the corollary.

\sect{Existence of dissipative solutions of ELHM}\label{ExistDiss}
In this section, we will use the notation below, as is done in \cite{Lions}:
\begin{align*}
E(v)\coloneq{}&-f_v \\
\norm{D^-v}_\8\coloneq{}&\norm{\pa{\sup_{|\jg|=1}\pa{-\xfr12(\pd_iv_j+\pd_jv_i)\jg,\jg}\2}^+}_{L^\8(D)}.
\end{align*}
The aim of this section is to prove the following.
\xbegin{teor}[Existence of dissipative solutions]
Let $u_0\in L^2(D)$ be divergence-free and $\log n_0\in L^{2+\eg}$ for some $\eg>0$. Then there exists at least one dissipative solution $u\in L^2(D_T;L^2_w(D)\!)$ of \eqref{eq:ELHM} with initial datum $u_0$.
\xend{teor}
\begin{qeddim}
We start, as in \cite{Lions}, from solutions to the viscous version, as provided by \kcref{thm:teor:LLR}:
\begin{equation}
\pd_t(u_{\ng}-\Dg^{-1}u_{\ng})+(u_{\ng}\per\grad)u_{\ng}-u_{\ng}^\perp\log n_0-\ng\Dg u_{\ng}=0, \label{eq:VELHMnu}
\end{equation}
with datum $\lbar u$ and the estimate:
\begin{equation}
\xfr12\Der{t}\xints D{}\abs{u_{\ng}}^2+\abs{(-\Dg)^{-\xfr12}u_{\ng}}^2\diff x+\ng\xints D{}\abs{\grad u_{\ng}}^2\diff x\leq0. \label{eq:EEnu}
\end{equation}
We have $u_{\ng}\in L^2_tH^1_x\cap L^\8_t(L^2\cap H^{-1})_x\cap\s{C}^0_t(L^2_w)_x$, as per \kcref{thm:cor:WCont}, and in addition $u_{\ng}(t)\to\lbar u$ in $L^2$ as $t\to0$.
\nipar Consider a divergence-free field $v\in\s{C}^0_t\s{C}^1_x$ -- or in fact a smooth one, since it is easy to prove that requiring \eqref{eq:ELHMdiss} for $v\in\Cinf$ implies it for $v\in\s{C}^0_t\s{C}^1_x$ -- and assume $v|_{\pd D}\equiv0$. Multiplying \eqref{eq:VELHMnu} by $v$ and using $\opn{div}u_{\ng}=0$ and the fact $u_{\ng},v$ vanish on $\pd D$, we obtain:
\begin{align}
\Der{t}\xints D{}u_{\ng}\per(v-\Dg^{-1}v)\diff x+{\hsp{1cm}}& \notag \\
\ng\xints D{}\grad u_{\ng}\per\grad v\diff x={}&\xints D{}u_{\ng}\per\pd_t(v-\Dg^{-1}v)+(u_{\ng}\per\grad)v\per u_{\ng}+u_{\ng}^\perp\log n_0\per v\diff x. \label{eq:A}
\end{align}
We already noted above that $a^\perp\per b=-a\per b^\perp$, allowing us to rewrite \eqref{eq:A} as:
\begin{equation}
\Der{t}\xints D{}u_{\ng}\per(v-\Dg^{-1}v)\diff x+\ng\xints D{}\grad u_{\ng}\per\grad v\diff x=\xints D{}u_{\ng}\per\pa{-E(v)+(u_{\ng}-v)\per\grad v}\diff x. \label{eq:Aa}
\end{equation}
In addition:
\begin{align}
\xfr12\Der{t}\xints D{}\abs v^2+\abs{(-\Dg)^{-\xfr12}v}^2\diff x={}&\xints D{}v\pd_t(v-\Dg^{-1}v)\diff x=\xints D{}v\sq{-E(v)-(v\per\grad)v+v^\perp\log n_0}\diff x={} \notag \\
{}={}&-\xints D{}v\per E(v)\diff x, \label{eq:B}.
\end{align}
provided $\opn{div}v=0$ (otherwise we have an additional $+\int_D|v|^2\opn{div}v\diff x$). Thus, combining \eqref{eq:EEnu}, \eqref{eq:Aa}, and \eqref{eq:B}, we conclude:
\begin{align}
\xfr12\Der{t}\xints D{}\abs{u_{\ng}-v}^2+\abs{(-\Dg)^{-\xfr12}(u_{\ng}-v)}^2\diff x\leq{}&\xints D{}u_{\ng}\per\pa{E(v)-(u_{\ng}-v)\per\grad v}\diff x+{} \notag \\
&{}-\ng\xints D{}\abs{\grad u_{\ng}}^2\diff x+\ng\xints D{}\grad u_{\ng}\per\grad v\diff x-\xints D{}v\per E(v)\diff x\leq{} \notag \\
{}\pux[\Big]{$\int_Dv\per[(u_{\ng}-v)\per\grad]v\diff x=0$ since $\opn{div}(u_{\ng}-v)=0$}{\leq}{}&\xints D{}E(v)(u_{\ng}-v)+\ng\grad u_{\ng}\per\grad v-[(u_{\ng}-v)\per\grad]v\per(u_{\ng}-v)\diff x\leq{} \notag \\
{}\leq{}&\norm{D^-v}_\8\xints D{}\abs{u_{\ng}-v}^2\diff x+\ng\norm{v}_{\dot H^1}\norm{u_{\ng}}_{\dot H^1}+{} \notag \\
&{}+\xints D{}E(v)(u_{\ng}-v)\diff x. \label{eq:D}
\end{align}
Integrating \eqref{eq:EEnu} in $\diff t$, we obtain:
\[\xfr12\xints D{}\abs{u_{\ng}}^2+\abs{(-\Dg)^{-\xfr12}u_{\ng}}^2\diff x\eval{T}+\ng\xints0T\xints D{}\abs{\grad u_{\ng}}^2\diff x\leq\xfr12\xints D{}\abs{\lbar u}^2+\abs{(-\Dg)^{-\xfr12}\lbar u}^2\diff x,\]
which implies, by concavity of $x\mapsto\rad x$ and ``reverse Jensen'':
\begin{equation}
\ng\xints0T\norm{u_{\ng}}_{\dot H^1}\diff x\leq\ng^{\xfr12}\pa{\ng\xints0T\xints D{}\abs{\grad u_{\ng}}^2\diff x\diff s}^{\xfr12}\leq C_T\ng^{\xfr12}, \label{eq:CTNu}
\end{equation}
where $C_T$ depends on $T$ but not $\ng$. Since $v$ is smooth and thus $\max\{\|D^-v\|_\8,\|v\|_{\dot H^1}\}\leq C$ and thus $\max\{\int_0^T\|D^-v\|_\8,\int_0^T\|v\|_{\dot H^1}\}\leq CT$, manipulating \eqref{eq:D} as was done a few sections ago with \eqref{eq:D1} and using the above estimate \eqref{eq:CTNu}, we conclude:
\begin{align}
\xints D{}\abs{u_{\ng}-v}^2+\abs{(-\Dg)^{-\xfr12}(u_{\ng}-v)}^2\diff x\leq{}&e^{2\int_0^t\|D^-v\|_\8\diff s}\xints D{}\abs{\lbar u-\lbar v}^2+\abs{(-\Dg)^{-\xfr12}(\lbar u-\lbar v)}^2\diff x+{} \notag\\
&{}+2\xints0T\xints D{}e^{2\int_s^t\|D^-v\|_\8\diff\tg}E(v)(u_{\ng}-v)\diff x\diff s+\tilde C_T\ng^{\xfr12}. \label{eq:DissELHM}
\end{align}
Provided, then, that $u_{\ng}$ converges in a way that allows one to pass to the limit for $\ng\to0$ in \eqref{eq:DissELHM}, the limit will be a dissipative solution of \eqref{eq:ELHM}. In particular, $u_{\ng}\xrightharpoonup{L^2}u$ is enough, as we will see below. 
\nipar To obtain the weak convergence, note that:
\[\pd_t(u_{\ng}-\Dg^{-1}u_{\ng})=-\pd_j(u_{\ng,j}u)+u_{\ng}^\perp\log n_0+\ng^{\xfr12}\pd_j(\ng^{\xfr12}\pd_ju_{\ng}).\]
The RHS is bounded in $L^2_tH^{-1}_x+L^\8_tW^{-(1+s),1}_x$ for all $s>0$. Indeed, term 1 is bounded in $L^\8_tW^{-(1+s),1}_x$, term 2 is bounded in $L^\8_tL^2_x\hra  L^2_tH^{-1}_x$ (since $\log n_0\in L^{2+\eg}$), and term 3 is bounded in $L^2_tH^{-1}_x$.
\nipar By reflexivity of $L^2$ and \cite[Lemma C.1]{Lions} (whose uniform continuity condition is implied by the above-stated property of $\pd_t(u_{\ng}-\Dg^{-1}u_{\ng})$), we conclude $u_{\ng}$ is precompact in $\s{C}^0_t(L^2_x)_w$, meaning a subsequence of $u_{\ng}-\Dg^{-1}u_{\ng}$ converges to some $w\in\s{C}^0_tL^2_x$ weakly uniformly in time.
\nipar $\Dg^{-1}$ is continuous $L^2\to H^2$, at least on the torus, thus compact $L^2\to L^2$. Since $u_{\ng}$ is bounded in $L^2$, $\Dg^{-1}u_{\ng}$ has a convergent subsequence. This means that, for a subsequence that satisfies both convergences, we must have $u_{\ng}\rightharpoonup U$ in $L^2$. Given the continuity, if $u_{\ng}\rightharpoonup U$, then $\Dg^{-1}u_{\ng}\to\Dg^{-1}U$.
\nipar We can then note that $\int|u-v|^2\lsim\liminf_n\int|u_{\ng}-v|^2$ by weak convergence, integrate by parts the $(-\Dg)^{-\xfr12}$ term and note that $u_{\ng}\rightharpoonup U$ and $\Dg^{-1}u_{\ng}\to\Dg^{-1}U$ in $L^2$, so that term passes to the limit, and weak convergence allows us to pass to the limit on the RHS. Thus, $U$ is a dissipative solution.
\end{qeddim}

\appendix
\sect{Some relations} \label{AppA}
\xbegin{lemma}
For any two vectors $a,b\in\R^2$, we have:
\begin{align*}
a\per b^\perp={}&-a\per b^\perp \\
a^\perp\per b^\perp={}&a\per b \\
(a^\perp)^\perp={}&-a.
\end{align*}
\xend{lemma}
\begin{qeddim}
If the latter two hold, then:
\[a^\perp\per b=(a^\perp)^\perp\per b^\perp=(-a)\per b^\perp,\]
so the first one holds.
\nipar $(a^\perp)^\perp=(-a_2,a_1)^\perp=(-a_1,-a_2)=-a$, proving the last relation, and $a^\perp\per b^\perp=(-a_2,a_1)\per(-b_2,b_1)=(-a_2)(-b_2)+a_1b_1=a_1b_1+a_2b_2=a\per b$, proving the first relation.
\end{qeddim}
\xbegin{lemma}
Let $\fg$ be a 2-variable scalar function and $u$ a 2-variable vector field. Define:
\begin{align*}
\fg_{3D}(x,y,z)\coloneq{}&\fg(x,y) \\
u_{3D,i}(x,y,z)\coloneq{}&u_i(x,y)\quad i=1,2 \\
u_{3D,3}(x,y,z)\coloneq{}&0 \\
\opn{curl}u\coloneq{}&\pd_1u_2-\pd_2u_1=\text{``}\grad^\perp\per u\text{''} \\
\grad^\perp\fg\coloneq{}&(-\pd_2\fg,\pd_1\fg) \\
(\opn{curl}\fg)_{3D}\coloneq{}&\opn{curl}(0,0,\fg_{3D})=(-\grad^\perp\fg)_{3D}.
\end{align*}
Then:
\begin{align}
\grad^\perp\fg={}&-\grad\fg_{3D}\x \hat z \label{eq:A1} \\
(\opn{curl}u)\hat z={}&\opn{curl}u_{3D} \label{eq:A2} \\
\opn{curl}\opn{curl}u={}&\grad\opn{div}u-\Dg u \label{eq:A3} \\
\opn{div}\grad^\perp\fg={}&0 \label{eq:A4} \\
\opn{curl}\grad\fg={}&0 \label{eq:A5} \\
\opn{div}u={}&\opn{curl}u^\perp \label{eq:A6} \\
\grad^\perp\opn{curl}u={}&\Dg u-\grad\opn{div}u. \label{eq:A7}
\end{align}
\xend{lemma}
\begin{qeddim}
There are all proved by simple calculations using the definitions.
\begin{align*}
\grad^\perp\fg={}&(-\pd_2\fg,\pd_1\fg)=(-\pd_2\fg_{3D},\pd_1\fg_{3D})=\det\pa{\mat{ccc} \hat x & \hat y & \hat z \\ -\pd_1\fg_{3D} & -\pd_2\fg_{3D} & 0 \\ 0 & 0 & 1 \emat}={} \\
{}={}&-\grad\fg_{3D}\x\hat z,
\intertext{which proves \eqref{eq:A1}.}
\opn{curl}u_{3D}={}&\det\pa{\mat{ccc} \hat x & \hat y & \hat z \\ \pd_1 & \pd_2 & \pd_3 \\ u_{3D,1} & u_{3D,2} & 0 \emat}={} \\
{}={}&(\pd_20-\pd_3u_{3D,2})\hat x-(\pd_10-\pd_3u_{3D,1})\hat y+(\pd_1u_{3D,2}-\pd_2u_{3D,1})\hat z={} \\
{}={}&(\pd_1u_2-\pd_2u_1)\hat z=(\opn{curl}u)\hat z, \\
\intertext{which proves \eqref{eq:A2}.}
(\opn{curl}\opn{curl}u)_{3D}={}&\opn{curl}\opn{curl}u_{3D}=\sum_{ijk\ell m}\eg_{ijk}\pd_j[\eg_{k\ell m}\pd_\ell u_{3D,m}]e_i={} \\
{}={}&\sum_{ij\ell m}(\dg_{i\ell}\dg_{jm}-\dg_{j\ell}\dg_{im})\pd_j\pd_\ell u_{3D,m}e_i={} \\
{}={}&\sum_{ij}(\pd_j\pd_iu_{3D,j}-\pd_j\pd_ju_{3D,i})e_i=\sum_i[\pd_i\opn{div}u_{3D}-\Dg u_{3D,i}]e_i{} \\
{}={}&\grad\opn{div}u_{3D}-\Dg u_{3D} \\
\intertext{which proves \eqref{eq:A3}. We can also see this in two other ways. The first one remains in 2D:}
\opn{curl}\opn{curl}u={}&\opn{curl}(\pd_1u_2-\pd_2u_1)=[\pd_2(\pd_1u_2-\pd_2u_1),-\pd_1(\pd_1u_2-\pd_2u_1)]={} \\
{}={}&[\pd_1\opn{div}u-\pd_1\pd_1u_1-\pd_2\pd_2u_1,\pd_2\opn{div}u-\pd_2\pd_2u_2-\pd_1\pd_1u_2]={} \\
{}={}&\grad\opn{div}u-\Dg u, \\
\intertext{The second other way is to derive it from a fully 3D relation:}
\opn{curl}\opn{curl}v={}&\opn{curl}\det\pa{\mat{ccc} \hat x & \hat y & \hat z \\ \pd_1 & \pd_2 & \pd_3 \\ v_1 & v_2 & v_3 \emat}=\opn{curl}\pa{\mat{c} \pd_2v_3-\pd_3v_2 \\ \pd_3v_1-\pd_1v_3 \\ \pd_1v_2-\pd_2v_1 \emat}={} \\
{}={}&\det\pa{\mat{ccc} \hat x & \hat y & \hat z \\ \pd_1 & \pd_2 & \pd_3 \\ \pd_2v_3-\pd_3v_2 & \pd_3v_1-\pd_1v_3 & \pd_1v_2-\pd_2v_1 \emat}={} \\
{}={}&\pa{\mat{c} \pd_2(\pd_1v_2-\pd_2v_1)-\pd_3(\pd_3v_1-\pd_1v_3) \\ \pd_3(\pd_2v_3-\pd_3v_2)-\pd_1(\pd_1v_2-\pd_2v_1) \\ \pd_1(\pd_3v_1-\pd_1v_3)-\pd_2(\pd_2v_3-\pd_3v_2) \emat}={} \\
{}={}&\pa{\mat{c} \pd_1(\opn{div}u-\pd_1v_1)-(\Dg v_1-\pd_1^2v_1) \\ \pd_2(\opn{div}u-\pd_2v_2)-(\Dg v_2-\pd_2^2v_2) \\ \pd_3(\opn{div}u-\pd_3v_3)-(\Dg v_3-\pd_1^3v_3) \emat}={} \\
{}={}&\grad\opn{div}u-\Dg u.
\intertext{}
\opn{div}\grad^\perp\fg={}&\opn{div}(-\pd_2\fg,\pd_1\fg)=-\pd_1\pd_2\fg+\pd_2\pd_1\fg=0 \\
\intertext{which proves \eqref{eq:A4}.}
\opn{curl}\grad\fg={}&\opn{curl}(\pd_1\fg,\pd_2\fg)=\pd_1\pd_2\fg-\pd_2\pd_1\fg=0 \\
\intertext{which proves \eqref{eq:A5}.}
\opn{curl}u^\perp={}&\opn{curl}(-u_2,u_1)=\pd_1u_1-\pd_2(-u_2)=\opn{div}u \\
\intertext{which proves \eqref{eq:A6}.}
\grad^\perp\opn{curl}b+\grad\opn{div}b={}&\pa{\mat c \cancel{-\pd_2\pd_1b_2}+\pd_2\pd_2b_1+\pd_1\pd_1b_1+\cancel{\pd_1\pd_2b_2} \\ \pd_1\pd_1b_2\bcancel{{}-\pd_1\pd_2b_1+\pd_2\pd_1b_1}+\pd_2\pd_2b_2 \emat}=\pa{\mat c\Dg b_1 \\ \Dg b_2 \emat}=\Dg b,
\end{align*}
which proves \eqref{eq:A7}.
\end{qeddim}
We used a summing rule for $\eg_{ijk}$, let us derive it. First and foremost:
\[\sum_k\eg_{ijk}\eg_{k\ell m}=\eg_{ij1}\eg_{1\ell m}+\eg_{ij2}\eg_{2\ell m}+\eg_{ij3}\eg_{3\ell m}.\]
If $i=j$ or $\ell=m$, all these terms are zero. If any of the four is the same as $k$, the corresponding term is zero. Thus, $\{i,j\}=\{\ell,m\}$, and $k\in\{1,2,3\}\ssm\{i,j\}$ is the only term that remains. This means the sum above is always either 1 or $-1$. Whichever $k$ may be, if $i=\ell,j=m$, then $\eg_{ijk}=\eg_{kij}=\eg_{k\ell m}$, thus the term is 1, and if $i=m,j=\ell$ then $\eg_{ijk}=-\eg_{kji}=-\eg_{k\ell m}$, so the term is $-1$. Now, if $i=\ell,j=m$, then $\dg_{i\ell}\dg_{jm}=1$ and $\dg_{im}\dg_{j\ell}=0$, otherwise the reverse holds. Thus, if $i=\ell,j=m$, then the sum above is $1=\dg_{i\ell}\dg_{jm}=\dg_{i\ell}\dg_{jm}-\dg_{im}\dg_{j\ell}$, whereas if $i=m,j=\ell$ then the sum is $-\dg_{j\ell}\dg_{im}$. Hence our rule.
\xbegin{lemma}[Integration by parts]
Let $a$ be a vector field and $b$ a scalar. Then:
\begin{align*}
\xints D{}b\opn{curl}a={}&-\xints D{}a\per\grad^\perp b\diff x+\xints{\mcl{\pd D}}{}(a\per\hg^\perp)b\diff\sg \\
\xints D{}\opn{curl}a\opn{curl}b\diff x={}&\xints D{}-(a\per\grad)\grad b+a\Dg b\diff x-\xints{\mcl{\pd D}}{}(a\per\hg^\perp)\grad^\perp b\diff\sg \\
\xints D{}\opn{curl}a\opn{curl}b\diff x={}&\xints D{}b\grad\opn{div}a-b\Dg a\diff x+\xints{\mcl{\pd D}}{}b\opn{curl}a\hg^\perp\diff\sg
\end{align*}
If instead $a,b$ are vector fields, then:
\[\xints D{}\opn{curl}a\per\opn{curl}b\diff x=-\xints D{}a\per[\Dg b-\grad\opn{div}b]\diff x+\xints{\mcl{\pd D}}{}(a\per\hg^\perp)\opn{curl}b\diff\sg.\]
\xend{lemma}
\begin{qeddim}
In all cases, we expand the curls and integrate by parts.
\begin{align*}
\xints D{}b\opn{curl}a\diff x={}&\xints D{}(\pd_1a_2-\pd_2a_1)b\diff x={} \\
{}={}&-\xints D{}a_2\pd_1b-a_1\pd_2b\diff x+\xints{\mcl{\pd D}}{}a_2b\hg_1-a_1b\hg_2\diff\sg={} \\
{}={}&\xints D{}a^\perp\per\grad b\diff x-\xints{\mcl{\pd D}}{}(a^\perp\per\hg)b\diff\sg={} \\
{}={}&-\xints D{}a\per\grad^\perp b\diff x+\xints{\mcl{\pd D}}{}(a\per\hg^\perp)b\diff\sg, \\
\intertext{which proves the first rule.}
\xints D{}\opn{curl}a\opn{curl}b\diff x={}&\xints D{}(\pd_1a_2-\pd_2a_1)\opn{curl}b\diff x={} \\
{}={}&-\xints D{}a_2\pd_1\opn{curl}b-a_1\pd_2\opn{curl}b\diff x+\xints{\mcl{\pd D}}{}(a\per\hg^\perp)\opn{curl}b\diff\sg={} \\
{}={}&\xints D{}a_2\pd_1\grad^\perp b-a_1\pd_2\grad^\perp b\diff x-\xints{\mcl{\pd D}}{}(a\per\hg^\perp)\grad^\perp b\diff\sg={} \\
{}={}&\xints D{}\pa{\mat{c} -a_2\pd_2\pd_1b+a_1\pd_2\pd_2b \\ a_2\pd_1\pd_1b-a_1\pd_1\pd_2b \emat}\diff x-\xints{\mcl{\pd D}}{}(a\per\hg^\perp)\grad^\perp b\diff\sg={} \\
{}={}&\xints D{}-(a\per\grad)\grad b+a\Dg b\diff x-\xints{\mcl{\pd D}}{}(a\per\hg^\perp)\grad^\perp b\diff\sg,
\intertext{which proves the second rule.}
\xints D{}\opn{curl}a\opn{curl}b\diff x={}&\xints D{}\opn{curl}a(\pd_1be_2-\pd_2be_1)\diff x={} \\
{}={}&-\xints D{}\pd_1\opn{curl}a\per be_2-\pd_2\opn{curl}a\per be_1\diff x+\xints{\mcl{\pd D}}{}b\opn{curl}a\hg^\perp\diff\sg={} \\
{}={}&-\xints D{}b\grad^\perp\opn{curl}a\diff x+\xints{\mcl{\pd D}}{}b\opn{curl}a\hg^\perp\diff\sg={} \\
{}={}&\xints D{}b\grad\opn{div}a-b\Dg a\diff x+\xints{\mcl{\pd D}}{}b\opn{curl}a\hg^\perp\diff\sg, \\
\intertext{which proves the third rule.}
\xints D{}\opn{curl}a\per\opn{curl}b\diff x={}&-\xints D{}a\per\grad^\perp\opn{curl}b\diff x+\xints{\mcl{\pd D}}{}(a\per\hg^\perp)\opn{curl}b\diff\sg={} \\
{}={}&-\xints D{}a\per[\Dg b-\grad\opn{div}b]+\xints{\mcl{\pd D}}{}(a\per\hg^\perp)\opn{curl}b\diff\sg,
\end{align*}
which proves the last rule. \placeqed*
\end{qeddim}
Note: for $a,b$ scalars, $\int_D\opn{curl}a\opn{curl}b\diff x=\int_D{}(-\grad^\perp a)\per(-\grad^\perp b)\diff x=\int_D\grad a\per\grad b\diff x$, so we just need to ``move a $\grad$'', which is a usual calculus rule.

\tableofcontents

\end{document}